\documentclass[a4paper, final, notitlepage, 10pt]{article}
\usepackage{latexsym,bm}
\usepackage{amsthm}
\usepackage{amsfonts}
\usepackage{amsmath}
\usepackage[shortalphabetic,abbrev]{amsrefs}
\usepackage[all]{xy}
\usepackage[pdftex,final]{graphics}
\usepackage[dvipsnames]{pstricks}
\usepackage{fancyhdr}
\usepackage{verbatim}
\usepackage{geometry}
\usepackage{xcolor}
\usepackage{enumerate}

\newtheorem{thm}{Theorem}[section]
\newtheorem{cor}[thm]{Corollary}
\newtheorem{lem}[thm]{Lemma}
\newtheorem{prop}[thm]{Proposition}

\theoremstyle{definition}

\newtheorem{rem}[thm]{Remark}
\newtheorem{exa}[thm]{Example}

\numberwithin{equation}{section}

\newcommand{\thmref}[1]{Theorem~\ref{#1}}

\newcommand{\lemref}[1]{Lemma~\ref{#1}}
\newcommand{\propref}[1]{Proposition~\ref{#1}}
\newcommand{\corref}[1]{Corollary~\ref{#1}}
\newcommand{\exaref}[1]{Example~\ref{#1}}
\newcommand{\remref}[1]{Remark~\ref{#1}}

\renewcommand{\O}{\mathcal{O}}

\newcommand{\B}{\overline{B}}

\newcommand{\N}{\overline{N}}
\renewcommand{\L}{\mathcal{L}}
\renewcommand{\o}{\Omega}
\renewcommand{\t}{\Theta}
\newcommand{\PP}{\mathbb{P}}

\newcommand{\bpi}{\overline{\pi}}
\newcommand{\e}{\varepsilon}
\newcommand{\nequiv}{\stackrel{num}{\sim}}
\newcommand{\tr}{\mathrm{tr}}

\title{Commuting involutions on surfaces of general type \\ with $p_g=0$ and $K^2=7$}
\author{Yifan Chen}
\date{}
\begin{document}
\maketitle
\renewcommand{\thefootnote}{\fnsymbol{footnote}}
\footnotetext{Yifan Chen: Academy of Mathematics and Systems Science, Chinese Academy of Sciences, No.~55 Zhongguancun East Road, Haidian District, Beijing, 100190, P.~R.~China}
\footnotetext{email address: chenyifan1984@gmail.com}
\begin{abstract}
 The aim of this article is to classify the pairs $(S, G)$, where $S$ is a smooth minimal surface of general type with $p_g=0$ and $K^2=7$, $G$ is a subgroup of the automorphism group of $S$
 and $G$ is isomorphic to the group $\mathbb{Z}_2^2$.
 The Inoue surfaces with $K^2=7$,
 which are finite Galois $\mathbb{Z}_2^2$-covers of the $4$-nodal cubic surface,
 are the first examples of such pairs.
 More recently, the author constructed  a new family of such pairs.
 They are finite Galois $\mathbb{Z}_2^2$-covers of certain $6$-nodal Del Pezzo surfaces of degree one.
 We prove that the base of the Kuranishi family of deformations of a surface in this family is smooth.
 We show that, in the Gieseker moduli space of canonical models of surfaces of general type,
 the subset corresponding to the surfaces in this family is an irreducible connected component,
 normal, unirational of dimension $3$.
\end{abstract}

\section{Introduction}
The first examples of surfaces of general type with $p_g=0$  are constructed in the 1930's
(cf.~\cite{Campedelli} and \cite{Godeaux}).
Since then, these surfaces have been studied by many mathematicians,
and more and more examples are constructed (cf.~\cite{BHPV}*{Table~14, page~304} and the references given there).
Nowadays, there is a long list of examples (cf.~\cite{survey}*{Tables 1-3}).
Minimal smooth surfaces of general type with $p_g=0$ have invariants $1 \le K^2 \le 9$.
Surprisingly, there are few examples of surfaces of general type with $p_g=0$ and $K^2=7$.
The first family of such surfaces is constructed by M.~Inoue (cf.~\cite{inoue}).

The bicanonical map plays an important role in the classification of surfaces of general type with $p_g=0$.
It is shown in \cite{bicanonical1} and \cite{bicanonical2} that the bicanonical morphism of a smooth minimal
surface of general type with $p_g=0$ and $K^2=7$ has degree $1$ or $2$; And if the bicanonical morphism has degree $2$,
the surface has a genus $3$ hyperelliptic fibration
and the fibration has five double fibers and one reducible fiber.
Another important way to classify surfaces with $p_g=0$ is
to study surfaces with certain automorphisms (for example, cf.~\cite{involution} and \cite{keumlee}).
Involutions on surfaces of general type with $p_g=0$ and $K^2=7$
are investigated in \cite{leeshin} and \cite{rito}.
All the possibilities of the quotient surfaces and the fixed loci of the involutions
are listed in these articles.

Recently, the author constructed a new family of surfaces of general type with $p_g=0$ and $K^2=7$ (cf.~\cite{chennew}).
In this article, we first demonstrate the process which leads to the discovery of these surfaces.
We explain the main idea.
Inspired by the results of \cite{leeshin} (See also~\cite{chennew}*{Section~6}),
we consider surfaces of general type with $p_g=0$, $K^2=7$ and with two distinct involutions.
We restrict our attention to the situation where the two involutions commute.

\begin{thm}\label{thm:classification}
Let $S$ be a minimal smooth surface of general type with $p_g=0$ and $K_S^2=7$.
Assume that $\mathrm{Aut}(S)$ contains a subgroup $G=\{1, g_1, g_2, g_3\}$,
which is isomorphic to $\mathbb{Z}_2^2$.
Let $R_i$ be the divisorial part of the fixed locus of the involution $g_i$ for $i=1, 2, 3$
and let $\pi \colon S \rightarrow \Sigma:=S/G$ be the quotient map.
Then the canonical divisor $K_S$ is ample and $R_i^2=-1$ for $i=1, 2, 3$.
Moreover, one of the following cases holds:
    \begin{enumerate}[\upshape (a)]
        \item $(K_SR_1, K_SR_2, K_SR_3)=(7, 5, 5)$; In this case,
              $S$ is an Inoue surface
              and $\pi \colon S \rightarrow \Sigma$ is a finite Galois $\mathbb{Z}_2^2$-cover of the
              $4$-nodal cubic surface as described in \exaref{exa:inoue};
        \item $(K_SR_1, K_SR_2, K_SR_3)=(5, 5, 3)$; In this case,
              the map $\pi \colon S \rightarrow \Sigma$ a finite Galois $\mathbb{Z}_2^2$-cover of a $6$-nodal Del Pezzo surface of degree one as described in \exaref{exa:new};
        \item $(K_SR_1, K_SR_2, K_SR_3)=(5, 3, 1)$;
              In this case, 
              $(R_1R_2, R_1R_3, R_2R_3)=(1, 3, 1)$ and
              the surface $\Sigma$ is a rational surface containing $8$ nodes and $K_{\Sigma}^2=-1$.
    \end{enumerate}
\end{thm}
We adopt the convention $K_SR_1 \ge K_SR_2 \ge K_SR_3$ in the theorem.

For simplicity, Inoue surfaces in the case (a) are referred to the surfaces with $p_g=0$ and $K^2=7$
constructed by M.~Inoue (\cite{inoue}).
The description of Inoue surfaces as $\mathbb{Z}_2^2$-covers of the $4$-nodal cubic surface appears in
\cite{bicanonical1}*{Example~4.1}.
It is also shown there that the bicanonical morphisms of Inoue surfaces have degree $2$.
The surfaces in the case (b) are constructed by the author (cf.~\cite{chennew}).
And they have birational bicanonical morphisms.

\thmref{thm:classification} is proved by \propref{prop:basicK2=7}, \thmref{thm:numerical}
and Section~4.
In particular,
Subsection~4.1 provides a detailed exposition of the classification process
which leads to the discovery of the surfaces in the case (b).
Here we explain the key strategy.
Because the Picard number of the surface $S$ is $3$,
the four divisors $K_S, R_1, R_2$ and $R_3$ are linearly dependent in the $\mathrm{Pic}(S)_{\mathbb{Q}}$.
Combining this fact, the algebraic index theorem and the adjunction formula,
we could easily analyze the configuration of the divisors $R_1, R_2, R_3$
and determine the number of nodes of the quotient surface $\Sigma$.
When $\Sigma$ is a Del Pezzo surface (the cases (a) and (b)), we get a complete classification.

However, we have difficulties when dealing with the case (c) because $\Sigma$ is no longer a Del Pezzo surface.
We do not know any examples for this case. We can not exclude this case at the moment.
Nevertheless, we make some remarks without proof.
In the case (c), from the classification table in \cite{leeshin},
one easily sees that the three intermediate quotient surfaces
$S/g_1, S/g_2$ and $S/g_3$ have Kodaira dimensions $0, 1$ and $2$ respectively.
In particular, the surface $S/g_3$ is a singular numerical Campedelli surface with five nodes.
Note that for surfaces in \exaref{exa:inoue} and \exaref{exa:new},
all the intermediate quotient surfaces have Kodaira dimensions at most $1$ (cf.~\cite{leeshin}*{Section~5}, \cite{chennew}*{Section~5 and Section~6}).
So it is worth finding pairs $(S, G)$ in the case (c).
We shall pursue this in the future.

The classification in \thmref{thm:classification} contributes to
the study of the moduli of the surfaces in the case (b).
The notation $\mathcal{M}^{\mathrm{can}}_{1, 7}$ in the following theorem stands for
the Gieseker moduli space of canonical models of surfaces of general type with $\chi(\O)=1$ and $K^2=7$
(cf.~\cite{Gieseker}).
\begin{thm}\label{thm:moduli}
Assume that $(S, G)$ is a pair in the case (b) of \thmref{thm:classification}.
Then the base of the Kuranishi family of deformations of $S$ is smooth.

In the Gieseker moduli space $\mathcal{M}^{\mathrm{can}}_{1,7}$,
the subset $\mathfrak{B}$ corresponding to the surfaces in the case (b) of \thmref{thm:classification} is an irreducible connected component, normal, unirational of dimension $3$.
\end{thm}
We point out that similar statements for Inoue surfaces have been achieved in \cite{inouemfd}:
The base of the Kuranishi family of deformations of an Inoue surface is smooth;
In the Gieseker moduli space,
the Inoue surfaces form a $4$-dimensional irreducible connected component,
normal  and unirational.
\thmref{thm:moduli} will be proved in Section~5.
The smoothness of the base of the Kuranishi family is obtained by calculating the dimensions of
the cohomology groups of the tangent sheaf of the surface $S$.
And the openness of $\mathfrak{B}$ in $\mathcal{M}^{\mathrm{can}}_{1, 7}$ follows from this.
We emphasize that the closedness of $\mathfrak{B}$ in $\mathcal{M}^{\mathrm{can}}_{1,7}$
follows from the classification in \thmref{thm:classification}.

\paragraph{\textbf{Notation and conventions}}
We adopt the convention that the indices $i \in \{1,2,3\}$ should be understood as residue classes modulo $3$.
We denote by $g_1,g_2,g_3$ the nonzero elements of the group $G \cong \mathbb{Z}_2^2$
and by $\chi_i \in G^*$ the nontrivial character orthogonal to $g_i$ for $i=1, 2, 3$.
Linear equivalence between divisors is denoted by $\equiv$.
Numerical equivalence between  divisors is denoted by $\nequiv$.
An $-m$-curve ($m \ge 0$) on a smooth projective surface
stands for a smooth rational curve with self intersection number $-m$.
A $-2$-curve is often called a nodal curve.
We denote by $c_1(\L)$ (respectively $c_1(D)$) the first Chern class of an invertible sheaf $\L$
(respectively a Cartier divisor $D$).
The rest of the notation is standard in algebraic geometry.

\section{Commuting involutions on surfaces with $p_g=0$}
Let $S$ be a smooth irreducible projective surface over $\mathbb{C}$.
A nontrivial automorphism $\alpha$ on $S$ is called an involution if $\alpha^2=\mathrm{Id}_S$.
We refer to \cite{involution}*{Section~3} for the properties of an involution on a surface.
We follow the ideas and the techniques there to study commuting involutions
on a surface of general type with $p_g=0$ and $K^2=7$.

\subsection{General case}To study the general case,
we first assume that $S$ is a smooth minimal surface of general type with $p_g=0$
and $\mathrm{Aut}(S)$ contains a subgroup $G=\{1, g_1, g_2, g_3\}$, which is isomorphic to $\mathbb{Z}_2^2$.
Burniat surfaces are examples of such surfaces with $2 \le K^2 \le 6$ (cf.~\cite{burniat} and \cite{peters}).

Let $\pi \colon S \rightarrow \Sigma = S /G$ be the quotient map.
We use Cartan's lemma (see~\cite{cartan}) to analyze the local properties of
the ramification locus and branch locus of $\pi$.
More precisely, for $i=1, 2, 3$,
let $R_i$ be the divisorial part of the fixed locus of the involution $g_i$ and let $B_i:=\pi(R_i)$.
Cartan's lemma implies that the divisors $R_1, R_2$ and $R_3$ satisfy the following properties:
\begin{enumerate}[\upshape (i)]
    \item if $R_i$ is not $0$, it is a disjoint union of irreducible smooth curves;
    \item the divisor $R_1+R_2+R_3$ is normal crossing.
\end{enumerate}
For each $i\in \{1,2,3\}$,  if $R_{i+1} \cap R_{i+2} \not =\emptyset$,
then the intersection points of $R_{i+1}$ and $R_{i+2}$ are isolated fixed points of the involution $g_i$.
The image of these points under $\pi$ are smooth points of $\Sigma$.
Besides the intersection points of $R_{i+1}$ and $R_{i+2}$,
the involution $g_i$ has $l_i$  pairs of isolated fixed points $(p_{ij_i},q_{ij_i})$ for $j_i=1, \ldots, l_i$.
The two points $p_{ij_i}$ and $q_{ij_i}$ of such a pair are permutated by $g_{i+1}$ and $g_{i+2}$.
Their images $r_{ij_i}=\pi(p_{ij_i})=\pi(q_{ij_i})$ are nodes of $\Sigma$.
The nodes $r_{ij_i}$ ($j_i=1, \ldots, l_i$ and $i=1, 2, 3$) are the only singularities of $\Sigma$.
In particular, $\Sigma$ is Gorenstein. We have the following formula
\begin{align}\label{eq:adjunction}
K_S=\pi^*K_{\Sigma}+R_1+R_2+R_3.
\end{align}

\begin{rem}\label{rem:B}
The discussion above also shows that the divisors $B_1, B_2$ and $B_3$ are contained in the smooth locus of $\Sigma$.
Moreover, the statements (i) and (ii) still hold if we replace $R_i$ by $B_i$
(cf.~\cite{singularbidouble}*{Theorem~2}).
\end{rem}

\begin{prop}\label{prop:basic1}
    Let $S$ be a minimal smooth surface of general type with $p_g=0$.
    Assume that $\mathrm{Aut}(S)$ contains a subgroup $G=\{1, g_1, g_2, g_3\}$,
    which is isomorphic to $\mathbb{Z}_2^2$. Then:
       \begin{enumerate}[\upshape (a)]
           \item for $i=1, 2, 3$, $2l_i+R_{i+1}.R_{i+2}=K_S.R_i+4$, $0 \le K_SR_i \le K_S^2$
                 and the integers $K_SR_i$, $R_{i+1}R_{i+2}$ and $K_S^2$ are of the same parity;
           \item $h^0(S, \O_S(2K_S))^{\mathrm{inv}}=\frac 14(K_S^2+K_SR_1+K_SR_2+K_SR_3)+1$ and for $i=1, 2, 3$,\\
                 $h^0(S, \O_S(2K_S))^{\mathrm{\chi_i}}=\frac{1}{4}(K_S^2+K_SR_i-K_SR_{i+1}-K_SR_{i+2})$;
           \item $K_\Sigma^2 =\frac 14 (K_S^2+R_1^2+R_2^2+R_3^2)+6-l_1-l_2-l_3$ and
                 the integer $K_S^2+R_1^2+R_2^2+R_3^2$ is divisible by $4$.
       \end{enumerate}
\end{prop}
\proof
The discussion above shows that the number of the isolated fixed points of
the involution $g_i$ is $2l_i+R_{i+1}R_{i+2}$.
By \cite{involution}*{Lemma~3.2 and Proposition~3.3~(v)},
$2l_i+R_{i+1}R_{i+2}=K_SR_i+4 \le K_S^2+4$
and $2l_i+R_{i+1}R_{i+2}$ has the same parity of $K_S^2$.
This implies (a).

Fix $i$.
Classical formulae for double cover show that the invariant subspace of $H^0(S, \O_S(2K_S))$ for the action of $g_i$
has dimension $\frac 12 (K_S^2+K_SR_i)+1$.
In our notation, we have
\begin{align*}
\dim H^0(S, \O_S(2K_S))^{\mathrm{inv}}+\dim H^0(S, \O_S(2K_S))^{\mathrm{\chi_i}}=\frac 12 (K_S^2+K_SR_i)+1\ \text{for}\ i=1, 2, 3.
\end{align*}
Since $\dim H^0(S, \O_S(2K_S))=K_S^2+1$,
assertion (b) follows.

By \eqref{eq:adjunction}, we have
\begin{align}
K_\Sigma^2=\frac 14 (K_S^2+R_1^2+R_2^2+R_3^2)-\frac 12 K_S(R_1+R_2+R_3)+\frac 12 (R_1R_2+R_1R_3+R_2R_3) \label{eq:KSigma2}
\end{align}
Because $\Sigma$ has only nodes, $K_\Sigma^2$ is an integer.
Then assertion (c) follows from (a).
\qed

\begin{cor}[\cite{involution}*{Corollary~3.4}]\label{cor:basic2}
    We keep the assumption and notation introduced above and fix $i \in \{1, 2 , 3\}$.
    The bicanonical map $\varphi \colon S \dashrightarrow \PP^{K_S^2}$ is composed with the involution $g_i$
    if and only if $K_SR_i=K_S^2$. In this case, $K_SR_{i+1}=K_SR_{i+2}$.
\end{cor}

Let $\eta \colon W \rightarrow \Sigma$ be the minimal resolution of $\Sigma$.
For each $i=1, 2, 3$,
let $\N_i$ be the disjoint union of the $l_i$ nodal curves over the nodes $r_{j_i}$ for $j_i=1,\ldots,l_i$
and let $\B_i:=\eta^*B_i$.
Let $\e \colon V \rightarrow S$ be the blowup at  $p_{j_i}$ and $q_{j_i}$ for $i=1, 2, 3$ and $j_i=1, \ldots, l_i$.
Then the $\mathbb{Z}_2^2$-action on $S$ lifts to $V$ and $V/G \cong W$.
There is a commutative diagram:
\begin{align}\label{diag}
\xymatrix@M=0em{
    &                                              &V \ar^{\bpi}"1,4" \ar_{\e}"2,3" & W \ar^{\eta}"2,4"
    & \B_i \ar_{\eta}"2,5" & {\quad}\N_i{\quad} \ar_{\eta}"2,6"\\
R_i & \{(p_{ij_i}, q_{ij_i})\}_{j_i=1, \ldots, l_i}  &  S \ar_{\pi} "2,4"             & \Sigma      & B_i=\pi(R_i)  & \{r_{j_i}=\pi(p_{ij_i})\}_{j_i=1, \ldots, l_i}
}
\end{align}
The quotient map $\bpi \colon V \rightarrow W$ is a finite flat $\mathbb{Z}_2^2$-cover branched on
the divisors $\Delta_1:=\B_1+\N_1, \Delta_2:=\B_2+\N_2$ and $\Delta:=\B_3+\N_3$.
There are three divisors $\L_1,\L_2$ and $\L_3$ of $W$ such that
\begin{align}
    2\L_i \equiv \B_{i+1}+\N_{i+1}+\B_{i+2}+\N_{i+2},\ \
    \L_i+\B_{i}+\N_{i} \equiv \L_{i+1}+\L_{i+2} \label{eq:buildingdata}
\end{align}
for $i=1, 2, 3$ (cf.~\cite{singularbidouble}*{Section~2}).

\begin{lem}\label{lem:DM}
Let
    \begin{align}\label{eq:DM}
          D:=2K_W+\B_1+\B_2+\B_3\ \text{and}\ M:=K_W+D
    \end{align}
Then:
    \begin{enumerate}[\upshape (a)]
        \item the divisor $D$ is nef and big, $\bpi^*D=\e^*(2K_S)$ and $D^2=K_S^2$;
        \item for each $i=1, 2, 3$, $\B_i^2=R_i^2$, $\B_i\B_{i+1}=R_iR_{i+1}$ and $D\B_i=K_SR_i$;
        \item $h^0(W, \O_W(D))=\frac 14 (K_S^2+K_SR_1+K_SR_2+K_SR_3)+1$ and $\dim|D| \ge 1$;
        \item any divisor in $|D|$ is $1$-connected and $p_a(D)\ge 1$;
        \item $h^0(W, \O_W(M))=p_a(D) \ge 1$.
    \end{enumerate}
    \noindent Assume furthermore that $K_S$ is ample, then:
    \begin{enumerate}[\upshape (a)]
    \item[\upshape (f)] if $DC=0$ for an irreducible curve $C$,
               then $C$ is one of the nodal curves in $\N_1 \cup \N_2 \cup \N_3$;
    \item[\upshape (g)] the divisor $M$ is nef.
    \end{enumerate}
\end{lem}
\proof Note that $\pi^*B_i=2R_i$ for $i=1, 2, 3$.
       By  \eqref{eq:adjunction} and from the diagram \eqref{diag}, we have
       $$\bpi^*D=\bpi^* \eta^* (2K_\Sigma+B_1+B_2+B_3)=\e^*\pi^*(2K_\Sigma+B_1+B_2+B_3)
       =\e^*(2K_S).$$
       Then $D^2=\frac 14 4K_S^2=K_S^2$. The divisor $D$ is nef and big because so is $2K_S$.

       Assertion~(b) follows from (a) and $\bpi^*\B_i=\e^*(2R_i)$.

       Note that
       \[K_W+\B_1+\B_2+\B_3 \nequiv \frac 12 D+ \frac 12 (\B_1+\B_2+\B_3).\]
       The divisor $D$ is nef and big, whereas $\frac 12(\B_1+\B_2+\B_3)$ is effective with zero integral part,
       and its support has normal crossings (see \remref{rem:B}).
       Thus $H^k(W, \O_W(D))=0$ for $k=1,2$ by Kawamata-Vieweg vanishing theorem (cf.~Corollary~5.12~(c) of \cite{ev}).
       The Riemann-Roch theorem gives $h^0(W, \O_W(D))=\frac 12 (D^2-D K_W)+1$.
       Then (c) follows from (b) and \eqref{eq:DM}.

        By \eqref{eq:DM}, (b) and \propref{prop:basic1}~(b),
\[\frac 12 (D^2+K_WD)=\frac 14 (3D^2-D\B_1-D\B_2-D\B_3)=\sum_{i=1}^{3}\dim H^0(S, \O_S(2K_S))^{\mathrm{\chi_i}}.\]
        So $\frac 12 (D^2 +K_WD) \ge 0$, i.e., $p_a(D) \ge 1$.
        Because $D$ is nef and big, $D$ is $1$-connected (cf.~Lemma~2.6 in \cite{adjoint}).
        This proves (d).

        Assertion~(e) follows from the long exact sequence obtained from
        \[0 \rightarrow \O_W(K_W) \rightarrow \O_W(K_W+D) \rightarrow \omega_{D} \rightarrow 0\]
        and the fact $p_g(W)=q(W)=0$.

        Now we prove (f).
        Assertion~(a) implies $\e^*K_S.\bpi^*C=0$ and thus $K_S. \e_{\ast}(\bpi^*C)=0$.
        Since $K_S$ is ample,
        $\e_{\ast}(\bpi^*C)=0$ and $\mathrm{Supp}\ \bpi^*C$ is contained in the exceptional divisors of $\e$.
        The nodal curves $\N_1 \cup \N_2 \cup \N_3$ are
        exactly the images of the exceptional divisors of $\e$ under $\bpi$.
        So $C$ is one of them.

        For (g), assume by contradiction that $MC<0$ for an irreducible curve $C$.
        Because $M$ is effective by (e), we have $C^2 <0$.
        Since $K_WC=(M-D)C<-DC \le 0$,
        $C$ is a $(-1)$-curve and $D.C=0$. This contradicts (f).
        Hence $M$ is nef.
\qed

\subsection{Surfaces with $p_g=0$ and $K^2=7$}
In the remainder of the article,
we always assume that $S$ is a smooth minimal surface of general type with $p_g=0$ and $K_S^2=7$.
We list some basic properties of $S$.
The surface $S$ has irregularity $q(S)=0$ and has Picard number $\rho(S)=3$ by Noether's formula and Hodge decomposition.
The expotential cohomology sequence gives $\mathrm{Pic}(S)\cong H^2(S, \mathbb{Z})$.
Poincar\'e duality implies that the intersection form on
$\mathrm{Num}(S):=\mathrm{Pic}(S)/\mathrm{Pic}(S)_{\mathrm{tor}}$ is unimodular.
The the bicanonical map $\varphi \colon S \rightarrow \mathbb{P}^{7}$
has degree either $1$ or $2$ (cf.~\cite{bicanonical1} and \cite{bicanonical2}).
We also need the following lemmas.
\begin{lem}[cf.~\cite{keum}*{Theorem~1.4~(1)~(f)}]\label{lem:onenodalcurve}
The surface $S$ contains at most one nodal curve.
\end{lem}
\proof
       Assume by contradiction that $S$ contains two nodal curves $C_1$ and $C_2$.
       Then $C_1C_2 \le 1$. The matrix of the intersection numbers of $K_S, C_1$ and $C_2$ has determinant $21$ or $28$,
       either of which is not a square integer.
       Since $\rho(S)=3$, this contradicts the fact that the intersection form on $\mathrm{Num}(S)$ is unimodular.
\qed

\begin{lem}(See the proof of \cite{inouebloch}*{Proposition~3.6})\label{lem:oneinvolution}
Let $\alpha$ be an involution on $S$ and let $R_\alpha$ be the divisorial part of the fixed locus of $\alpha$.
Then $R_\alpha^2= \pm 1$.
\end{lem}
\proof Let $\tr(\alpha^*)$  be  the trace of the induced linear map
$\alpha^* \colon H^2(S, \mathbb{C}) \rightarrow H^2(S,\mathbb{C})$.
Then $R_\alpha^2=2-\tr(\alpha^*)$ by \cite{manynodes}*{Lemma~4.2}.
The Chern classes $c_1(K_S)$ and $c_1(R_\alpha) \in H^2(S, \mathbb{C})$ are invariant under $\alpha^*$.
So $\tr(\alpha^*)=-1, 1$ or $3$.
It suffices to exclude the case $\tr(\alpha^*)=-1$.

If $\tr(\alpha^*)=-1$, then $\alpha^*$ has eigenvalues $-1, -1$ and $1$.
This implies $R_\alpha \nequiv rK_S$ for some positive rational number $r$.
Because $R^2=2-\tr(\alpha^*)=3$ and $K_S^2=7$,
this is impossible.
\qed
\\

We assume furthermore that $\mathrm{Aut}(S)$ contains a subgroup $G=\{1, g_1, g_2, g_3\} \cong \mathbb{Z}_2^2$.
We keep the same notation introduced in the previous subsection and denote by
\begin{align}\label{eq:A}
A=(R_iR_j)_{1 \le i \le j \le 3}
\end{align}
the matrix of intersection numbers of the ramification divisors $R_1, R_2$ and $R_3$
of the quotient map $\pi \colon S \rightarrow \Sigma=S/G$.

Without loss of generality, we may assume $K_SR_1 \ge K_SR_2 \ge K_SR_3$.
Since the bicanonical map $\varphi$ has degree at most $2$,
$\varphi$ is composed with at most one involution in $G$.
Then by \propref{prop:basic1}~(a), (b) and \corref{cor:basic2},
one of the following cases occurs:
\begin{itemize}
\item if $\varphi$ is composed with exact one involution in $G$, then
      \begin{align}
      (K_SR_1, K_SR_2, K_SR_3) \in \{(7, 1, 1), (7, 3, 3), (7, 5, 5)\}; \label{eq:numerical1}
      \end{align}
\item if $\varphi_{2K_S}$ is not composed with any involution in $G$, then
      \begin{align}
      (K_SR_1, K_SR_2, K_SR_3) \in \{(3, 1, 1), (3, 3, 3), (5, 3, 1), (5, 5, 3)\}. \label{eq:numerical2}
      \end{align}
\end{itemize}

\begin{prop}\label{prop:basicK2=7}
    Let $S$ be a minimal smooth surface of general type with $p_g=0$ and $K_S^2=7$.
    Assume that $\mathrm{Aut}(S)$ contains a subgroup $G=\{1, g_1, g_2, g_3\}$.
    Then:
    \begin{enumerate}[\upshape (a)]
        \item for $i=1, 2, 3$, $R_i^2=-1$;
        \item the quotient $\Sigma$ is a rational surface with $K_\Sigma^2=7-l_1-l_2-l_3$ and $\rho(\Sigma)=3$;
        \item the divisors $R_1, R_2, R_3$ generate a sublattice of $\mathrm{Num}(S)$
              and $\det A$ is a positive square integer;
        \item the numbers $l_1, l_2$ and $l_3$ are even integers;
        \item the canonical divisor $K_S$ is ample.
    \end{enumerate}
\end{prop}
\proof
We remark that $R_i^2=\pm 1$ for $i=1, 2, 3$ by \lemref{lem:oneinvolution}.
We recall in the diagram \eqref{diag} that $W$ is the minimal resolution of $\Sigma$.
Note that $R_iR_{i+1} (=\B_i\B_{i+1})$ is a positive odd integer by
(see \propref{prop:basic1}~(a) and \lemref{lem:DM}~(b)).

We first exclude the case $(K_SR_1, K_SR_2, K_SR_3)=(3, 1, 1)$.
The algebraic index theorem and $K_S^2=7$ gives $(R_1+R_2+R_3)^2 \le \frac{5^2}{7}$.
Because $R_i^2=\pm 1$,  \propref{prop:basic1}~(a) implies $R_i^2=-1$ and $R_iR_{i+1}=1$ for $i=1, 2, 3$.
Then $\B_i^2=-1$, $\B_i\B_{i+1}=1$ and $D\B_2=1$ by \lemref{lem:DM}~(b).
So $K_W B_2=\frac 12(D-\B_1-\B_2-\B_3)\B_2=0$ by \eqref{eq:DM}.
This contradicts the adjunction formula.

Now we show that $W$ (and thus $\Sigma$) is a rational surface and $K_W^2 \le 3$.
By \eqref{eq:DM} and \lemref{lem:DM}~(b), we have
\begin{align}\label{eq:DKSigma}
D K_W=\frac 12 (D^2-D\B_1-D\B_2-D\B_3)=\frac 12 (K_S^2-K_SR_1-K_SR_2-K_SR_3)
\end{align}
It follows that $D K_W \in \{-5, -3, -1\}$
for all the possibilities \eqref{eq:numerical1} and \eqref{eq:numerical2}.
Since $D$ is nef and big (see \lemref{lem:DM}~(a)), $W$ is a rational surface.
The algebraic index theorem  and $D^2=7$ gives $K_W^2 \le 3$.

Because $W$ is a rational surface, the Picard number $\rho(W)=10-K_W^2 \ge 7$.
By \propref{prop:basic1}~(c),
one of the following two cases holds:
\begin{itemize}
    \item for $i=1, 2, 3$, $R_i^2=-1$ and $K_W^2=7-l_1-l_2-l_3$;
    \item two of the three integers $R_i^2$ equal $1$, the third one equals $-1$ and $K_W^2=8-l_1-l_2-l_3$.
\end{itemize}
Assume that the latter holds. Without loss of generality,
assume that $R_1^2=R_2^2=1$ and $R_3^2=-1$.
Then $\rho(W)=l_1+l_2+l_3+2$ and $W$ contains $\rho(W)-2$ disjoint nodal curves $\N_1 \cup \N_2 \cup \N_3$.
By \cite{manynodes}*{Theorem~3.3}, $l_1+l_2+l_3$ is an even integer.
Also $\B_1, \B_3$ and these nodal curves generate a sublattice of $\mathrm{Num}(W)$.
The matrix of intersection numbers of this sublattice has determinant $-2^{l_1+l_2+l_3}[1+(\B_1\B_3)^2]$.
Since the intersection form on $\mathrm{Num}(W)$ is unimodular,
$2^{l_1+l_2+l_3}[1+(\B_1\B_3)^2]$ is a  positive square integer.
Because $l_1+l_2+l_3$ is an even integer and $\B_1\B_3$ is a positive odd integer, this is impossible.
Thus the former case holds.
We have proved (a) and (b).

The matrix $A$ has determinant
\begin{align}\label{eq:detA}
    \det A = (R_1R_2)^2+(R_1R_3)^2+(R_2R_3)^2+2(R_1R_2)(R_1R_3)(R_2R_3)-1
\end{align}
Because the numbers $R_iR_{i+1}$ are positive odd integer, $\det A$ is a positive integer.
Since $\rho(S)=3$, $R_1, R_2, R_3$ generate a sublattice of $\mathrm{Num}(S)$ and
therefore $\det A$ is a square integer
by Poincar\'e duality.
This proves (c).

By \lemref{lem:DM}~(b), the matrix $A$ is also the intersection number matrix of the divisors $\B_1, \B_2, \B_3$.
It follows that $\B_1,\B_2$ and $\B_3$ and the $l_1+l_2+l_3$ nodal curves $\N_1 \cup \N_2 \cup \N_3$
generate a sublattice of $\mathrm{Num}(W)$.
Therefore $2^{l_1+l_2+l_3}\det A$ is a positive square integer by Poincar\'e duality.
Hence  $l_1+l_2+l_3$ is an even integer by (c).

Note that $R_i$ ($B_i$) is  a disjoint union of smooth irreducible curves (see \remref{rem:B}).
We apply the Hurwitz formula
for the double cover $\pi|_{R_i} \colon R_i \rightarrow B_i$ induced by the action of $g_{i+1}$ ($g_{i+2}$):
\[K_SR_i+R_i^2=2(2p_a(B_i)-2)+R_i(R_{i+1}+R_{i+2})\]
By \propref{prop:basic1}~(a), we have
\begin{align}\label{eq:genusB}
    2p_a(B_i)-2=\frac 12(K_SR_i+R_i^2-K_SR_{i+1}-K_SR_{i+2})+ l_{i+1}+l_{i+2}-4
\end{align}
For all the possibilities \eqref{eq:numerical1} and \eqref{eq:numerical2},
$K_SR_i+R_i^2-K_SR_{i+1}-K_SR_{i+2}$ is divisible by $4$ for each $i=1, 2, 3$.
So $l_{i+1}+l_{i+2}$ is an even integer for $i=1, 2, 3$ by \eqref{eq:genusB}.
We have seen that $l_1+l_2+l_3$ is an even integer.
Therefore $l_1, l_2, l_3$ are even integers.
This proves (d).

Now we prove (e).
Assume by contradiction that $K_S$ is not ample.
Then $S$ contains exactly one nodal curve $C$ by \lemref{lem:onenodalcurve}.
Thus $C$ is $G$-invariant.
Let $C'=\pi(C) \subset \Sigma$ and let $\bar{C'}$ be the strict transform of $C'$ on $W$.

First assume that $C$ is contained in some $R_i$, i.e.,
$C'$ is contained in some $B_i$.
By \remref{rem:B}, $C'$ is contained in the smooth locus of $\Sigma$ and $\pi^*C'=2C$.
Thus $C'^2=-2$.
It follows that $\bar{C'}$ is a nodal curve,
which is disjoint from the nodal curves $\N_1 \cup \N_2 \cup \N_3$.
Then $W$ contains $l_1+l_2+l_3+1$ pairwise disjoint nodal curves.
By (d), this contradicts \cite{manynodes}*{Theorem~3.3}.

Hence $C$ is not contained in $R_i$, i.e., $\B_i \not \ge \bar{C'}$ for $i=1, 2, 3$.
By \propref{lem:DM}~(a), $D\bar{C'}=0$ and $\bar{C'}^2 <0$.
It follows that $2K_W\bar{C'}=-(\B_1+\B_2+\B_3)\bar{C'} \le 0$.
Thus $\bar{C'}$ is either a $(-1)$-curve or a nodal curve.

If $\bar{C}'$ is nodal curve, then $\B_i\bar{C'}=0$ for $i=1, 2, 3$,
i.e., $R_iC=0$ of $i=1, 2, 3$.
This contradicts (c).
So $\bar{C}'$ is a $(-1)$-curve and $(\B_1+\B_2+\B_3)\bar{C'}=2$, i.e., $C'(B_1+B_2+B_3)=2$.
We remark that for any Galois $\mathbb{Z}_2^2$-cover $\mathbb{P}^1 \rightarrow \mathbb{P}^1$,
in the target space $\mathbb{P}^1$,
the branch locus consists of three distinct points.
Applying this remark to the cover $\pi|_C \colon C \rightarrow C'$,
because $C'$ intersects the divisorial part $B_1+B_2+B_3$ of the branch locus of $\pi$ at two points,
we conclude that $C'$ passes through exactly one node of $\Sigma$.
Equivalently,
$\bar{C'}$ intersects exactly one nodal curve $C_1$ in $\N_1 \cup \N_2 \cup \N_3$ and $\bar{C'}C_1=1$.
Blowing down $\bar{C'}$ and then the image of $C_1$, we obtain a surface $W'$ containing $l_1+l_2+l_3-1$
disjoint nodal curves and $\rho(W')=l_1+l_2+l_3+1$.
By (d), this again contradicts \cite{manynodes}*{Theorem~3.3}.

Hence $K_S$ is ample.
\qed

\begin{rem}\label{rem:BBBNNN}
      We see in the proof of \propref{prop:basicK2=7}~(c) that the three curves $\B_1,\B_2,\B_3$ and
      the $l_1+l_2+l_3$ nodal curves $\N_1\cup \N_2 \cup \N_3$ generate a sublattice of $\mathrm{Num}(W)$.
\end{rem}

The following theorem is the main result of this section.
\begin{thm}\label{thm:numerical}
    Let $S$ be a minimal smooth surface of general type with $p_g=0$ and $K_S^2=7$.
    Assume that $\mathrm{Aut}(S)$ contains a subgroup $G=\{1, g_1, g_2, g_3\}$.
    Let $\Sigma$ be the quotient surface $S/G$.
     Then one of the following cases occurs:
    \begin{enumerate}[\upshape (a)]
        \item $(K_SR_1, K_SR_2, K_SR_3)=(7, 5, 5)$;
              in this case, $(R_1R_2, R_1R_3, R_2R_3)=(5,9,7)$,  $(l_1, l_2, l_3)=(2, 0, 2)$
              and $K_\Sigma^2=3$;
        \item $(K_SR_1, K_SR_2, K_SR_3)=(5, 5, 3)$;
              in this case, one of the following cases occurs:
              \begin{enumerate}
                \item[(b1)] $(R_1R_2, R_1R_3, R_2R_3)=(7, 5, 1)$, $(l_1, l_2, l_3)=(4, 2, 0)$ and $K_\Sigma^2=1$;
                \item[(b2)] $(R_1R_2, R_1R_3, R_2R_3) \in \{(3, 5, 1), (7, 1, 1)\}$,
                            $l_1+l_2+l_3=8$ and $K_\Sigma^2=-1$;
              \end{enumerate}
        \item $(K_SR_1, K_SR_2, K_SR_3)=(5, 3, 1)$;
              in this case, $(R_1R_2, R_1R_3, R_2R_3)=(1, 3, 1)$, $(l_1, l_2, l_3)=(4, 2, 2)$ and $K_\Sigma^2=-1$.
    \end{enumerate}
\end{thm}
Here we keep the same notation introduced above and adopt the convention $K_SR_1 \ge K_SR_2 \ge K_SR_3$.
Furthermore, since $K_SR_2=K_SR_3$ in the case (a), we may assume $R_1R_2 \le R_1R_3$;
Similarly, we assume $R_1R_3 \ge R_2R_3$ in the case (b).

\proof All the possibilities are listed in \eqref{eq:numerical1}-\eqref{eq:numerical2}
       and the case $(3, 1, 1)$ has been excluded in the proof of \propref{prop:basicK2=7}.
       For each possibility, once we determine the matrix $A=(R_iR_j)_{1 \le i \le j \le 3}$,
       we can calculate $l_1, l_2, l_3$ and $K_\Sigma^2$
       by \propref{prop:basic1}~(a) and by \propref{prop:basicK2=7}~(b).
       Note that $R_i^2=-1$ by \propref{prop:basicK2=7}~(a).
       So it suffices to determine the intersection numbers $R_iR_{i+1}$ for $i=1, 2, 3$.

       First assume $(K_SR_1, K_SR_2, K_SR_3)=(7, 3, 3)$.
       The algebraic index theorem gives
                       \[2R_1R_2-2=(R_1+R_2)^2 \le \frac {(K_SR_1+K_SR_2)^2}{K_S^2}=\frac {10^2}{7}\]
       i.e., $R_1R_2 \le 8$. Similarly, we have
       $R_1R_3 \le 8 $, $R_2R_3 \le 3$ and $R_1R_2+R_1R_3+R_2R_3 \le 13$.
       By \propref{prop:basicK2=7}~(d) and \propref{prop:basic1}~(a), $K_SR_i-R_{i+1}R_{i+2}$ is divisible by $4$.
       Thus $R_1R_2, R_1R_3 \in \{3, 7\}$ and  $R_2R_3=3$.
       Then $\det A=80$ or $192$ by \eqref{eq:detA}.
       This contradicts \propref{prop:basicK2=7}~(c).
       Hence the case $(7, 3, 3)$ is excluded.

       The same reasoning will exclude the cases $(7, 1, 1)$ and $(3, 3, 3)$
       and gives $(R_1R_2, R_1R_3, R_2R_3)=(1, 3, 1)$ for the case $(5, 3, 1)$.
       This proves (c).

       For the case $(7, 5, 5)$, the algebraic index theorem gives $R_1R_2 \le 12$, $R_1R_3 \le 12$
       and $R_2R_3 \le 8$. Then \propref{prop:basicK2=7}~(d) and \propref{prop:basic1}~(a) imply
       $R_1R_2, R_1R_3 \in \{1, 5, 9\}$ and $R_2R_3 \in \{3, 7\}$.
       We shall use the \lemref{lem:DM}~(g) to narrow down the possibilities.
       On the minimal resolution $W$ of the quotient surface $\Sigma$,
       by \eqref{eq:DM} and \lemref{lem:DM}~(a)-(b), we have
       $DK_W=\frac 12D(D-\B_1-\B_2-\B_3)=-5$ and thus $M^2=(K_W+D)^2=K_W^2-3$.
       The algebraic index theorem yields $K_W^2 \le 3$.
       Because $K_S$ is ample, $M$ is nef (cf.~\lemref{lem:DM}~(g)).
       Then $M^2 \ge 0$ and thus $K_W^2 \ge 3$.
       Therefore $K_\Sigma^2=K_W^2=3$.
        We have $R_1R_2+R_1R_3+R_2R_3=21$ by \eqref{eq:KSigma2}.
       Then $(R_1R_2, R_1R_3, R_2R_3) \in \{(9, 9, 3), (5, 9, 7), (9, 5, 7)\}$.
       If $(R_1R_2, R_1R_3, R_2R_3)=(9, 9, 3)$, $\det A=656$ by \eqref{eq:detA},
       a contradiction to \propref{prop:basicK2=7}~(c).
       Hence $(R_1R_2, R_1R_3, R_2R_3)=(5,9,7)$ or $(9, 5, 7)$.
       This proves (a).

       For the case $(5, 5, 3)$,
       the algebraic index theorem, \propref{prop:basicK2=7}~(d) and \propref{prop:basic1}~(a)
       imply $R_1R_2 \in \{3, 7\}$ and $R_1R_3, R_2R_3 \in \{1, 5\}$.
       By \eqref{eq:DM} and \lemref{lem:DM}~(a)-(b),
       we have $DK_W=\frac 12(D-\B_1-\B_2-\B_3)=-3$ and $M^2=(K_W+D)^2=K_W^2+1$.
       The algebraic index theorem yields $K_W^2 \le 1$.
       Because $M$ is nef,  $M^2 \ge 0$ and thus $K_W^2 \ge -1$.
       Finally by \propref{prop:basicK2=7}~(b) and (d), $K_W^2= \pm 1$.

       If $K_W^2=1$, then $R_1R_2+R_1R_3+R_2R_3=13$ by \eqref{eq:KSigma2}.
       Then $(R_1R_2, R_1R_3, R_2R_3) \in \{(3, 5, 5), (7, 5, 1), (7, 1, 5)\}$.
       If $(R_1R_2, R_1R_3, R_2R_3)=(3, 5, 5)$, then  $\det A=208$ by \eqref{eq:detA},
       a contradiction to \propref{prop:basicK2=7}~(c).
       Hence $(R_1R_2, R_1R_3, R_2R_3)=(7, 5, 1)$ or $(7, 1, 5)$.

       If $K_W^2=-1$, then by \eqref{eq:KSigma2}, $R_1R_2+R_1R_3+R_2R_3=9$.
       Then $(R_1R_2, R_1R_3, R_2R_3) \in \{(3, 5, 1), (3, 1, 5), (7, 1, 1)\}$.

We complete the proof of \thmref{thm:numerical}.\qed\\

To conclude this section, we state a lemma on fibrations on $W$.
\begin{lem}\label{lem:xiao}
Assume that $|F|$ is a base point free pencil of  curves on $W$.
Then $DF \ge 2$.
\end{lem}
\proof
Let $F'=\e_{\ast}\bpi^*F$.
Since $\bpi^*D=\e^*2K_S$, it follows that $K_SF'=2DF$.
If $DF= 1$, then $K_SF'=2$.
And then $F'^2=0$ by the algebraic index theorem.
It follows that $|F'|$ induces a base point free pencil of genus $2$ on $S$.
This contradicts \cite{xiao}*{Theorem~2} because $K_S^2=7$.
\qed

\section{Examples}
We describe the known examples of surfaces of general type with $p_g=0$, $K^2=7$ and with commuting involutions.

\begin{exa}[Inoue surfaces as bidouble covers of the $4$-nodal cubic surface]\label{exa:inoue}
Inoue surfaces are the very first examples of surfaces of general type with $p_g=0$ and $K^2=7$ (cf.~\cite{inoue}).
The following description is from \cite{bicanonical1}
(see \cite{rito} for an equivalent description; cf.~\cite{chennew}*{Section~6}).
We change the notation for being coherent with Section~2.

Let $p_1,p_2,p_3,p_1',p_2',p_3'$ be the six vertices of a complete quadrilateral on $\PP^2$.
Let $\sigma \colon W \rightarrow \mathbb{P}^2$ be the blowup of these points.
Denote by $E_i$ (respectively $E_i'$) the exceptional curve of $W$ over $p_i$ (respectively $p_i'$)
and by $L$ the pullback of a general line by $\sigma$.
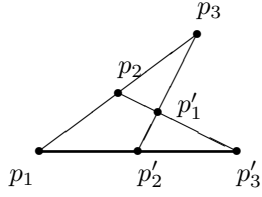
\begin{figure}[ht]
\setlength{\unitlength}{1.3mm}
\centerline{\begin{picture}(25,20)
\linethickness{1pt}
\put(0,0){$p_1$}   \put(3,3){\circle*{0.8}}
\put(3,3){\line(4,3){16}}
\put(3,3){\line(1,0){20}}
\put(11,11){$p_2$}  \put(11,9){\circle*{0.8}}
\put(11,9){\line(2,-1){12}}
\put(19,17){$p_3$}  \put(19,15){\circle*{0.8}}
\put(19,15){\line(-1,-2){6}}
\put(17,7){$p_1'$}  \put(15,7){\circle*{0.8}}
\put(13,0){$p_2'$} \put(13,3){\circle*{0.8}}
\put(23,0){$p_3'$} \put(23,3){\circle*{0.8}}
\end{picture}}
\caption{\footnotesize{Configurations of $p_1, p_1,\ldots, p_3'$} \label{figure1}}
\end{figure}

The surface $W$ has four disjoint nodal curves, their divisor classes are
\begin{align}
Z_i \equiv L-E_i-E_{i+1}'-E_{i+2}',\  Z \equiv L-E_1-E_2-E_3.
\end{align}
Note that $Z_1, Z_2, Z_3$ and $Z$ are the proper transforms of the four sides of the quadrilateral.
Let $\eta \colon W \rightarrow \Sigma$ be the morphism contracting there curves.
The surface $\Sigma$ has four nodes, $-K_\Sigma$ is ample and $\Sigma$ is the $4$-nodal cubic surface.

Let $\Gamma_1,\Gamma_2$ and $\Gamma_3$ be the proper transforms of the three diagonals of the quadrilateral,
i.e., $\Gamma_i \equiv L-E_i-E_i'$ for $i=1, 2, 3$.
For each $i=1,2,3,$
$W$ has a pencil of rational curves $|F_i|:=|2L-E_{i+1}-E_{i+2}-E_{i+1}'-E_{i+2}'|$.
Note that $-K_W \equiv \Gamma_1+\Gamma_2+\Gamma_3 \equiv \Gamma_i+F_i$ for $i=1,2,3$.

Now we define three effective divisors on $W$
\begin{align}
    \Delta_1:&=\Gamma_1+F_2+Z_1+Z_3,\nonumber\\
    \Delta_2:&=\Gamma_2+F_3,   \label{eq:inouedata}\\
    \Delta_3:&=\Gamma_3+F_1+F_1'+Z_2+Z. \nonumber
\end{align}
Here we require that $F_i$ ($i=1, 2, 3$) and $F_1'$ are smooth $0$-curves and the divisor
$\Delta: =\Delta_1+\Delta_2+\Delta_3$ has only nodes.
There is a smooth finite $G$-cover $\bpi \colon V \rightarrow W$
branched on the divisors $\Delta_1, \Delta_2$ and $ \Delta_3$.
The (set theoretic) inverse image $\bpi^{-1}Z_i$ or $\bpi^{-1}Z$ is a disjoint union of two $(-1)$-curves.
Let $\e \colon V  \rightarrow S$ be the blowdown of these eight $(-1)$-curves.
From the construction, there is a finite $G$-cover $\pi \colon S \rightarrow \Sigma$ such that
the diagram \eqref{diag} commutes.
And $S$ is a smooth minimal surface of general type with $p_g(S)=0$ and $K_S^2=7$.
Moreover, $K_S$ is ample and the bicanonical map of $S$ has degree two.
\end{exa}

In the notation of Section~2, $\B_1=\Gamma_1+F_2, \B_2=\Gamma_2+F_3, \B_3=\Gamma_3+F_1+F_1', \N_1=Z_1+Z_3, \N_2=0,
\N_3=Z_2+Z$ and $D=-K_W+F_1'$.
It follow that  $(D\B_1, D\B_2, D\B_3)=(7, 5, 5)$ and $(\B_1\B_2, \B_1\B_3, \B_2\B_3)=(5, 9, 7)$.
So the Inoue surfaces satisfy \thmref{thm:numerical}~(a) by \lemref{lem:DM}~(b).

\begin{rem}\label{rem:inouerem}
We will need the following remarks in the proof of \thmref{thm:classification}~(a). See Subsection~4.3.
\begin{enumerate}[\upshape (1)]
\item Note that $W$ contains exactly nine $(-1)$-curves $E_i, E_i'$ and $\Gamma_i$ for $i=1, 2, 3$.
      There are exactly three $(-1)$-curves $\Gamma_1, \Gamma_2$ and $\Gamma_3$,
      which are disjoint from the nodal curves.
\item If we replace $Z_1+Z_3$ in $\Delta_1$ by $Z_2+Z$ and replace $Z_2+Z$ in $\Delta_3$ by $Z_1+Z_3$,
      these new $\Delta_1, \Delta_2$ and $\Delta_3$ still
      define a smooth finite $\mathbb{Z}_2^2$-cover.
      In this way, we get another $4$-dimensional family of surfaces of general type.
      However, this family is the same as the original one.

      Indeed, let $\alpha$ be the involution on $\mathbb{P}^2$ such that
      $\alpha(p_k)=p_k'$ and $\alpha(p_k')=p_k$ for $k=1, 2$.
      Then $\alpha(p_3)=p_3$ and $\alpha(p_3')=p_3'$.
      It induces an involution $\alpha'$ on $W$.
      The $(-1)$-curves $\Gamma_1, \Gamma_2, \Gamma_3$ are $\alpha'$-invariant
      and the divisors classes of $F_1, F_2, F_3$ are $\alpha'$-invariant.
      We also have $\alpha'(Z_1)=Z_2, \alpha'(Z_2)=Z_1, \alpha'(Z_3)=Z$ and $\alpha'(Z)=Z_3$.

\item Observe that the two nodal curves in the same $\Delta_k$ ($k=1, 3$)
      are in the same singular member of the pencil $|F_2|$.
      Indeed, the singular members of the pencil $|F_2|$ are $\Gamma_1+\Gamma_3$, $Z_1+2E_2'+Z_3$
      and $Z_2+2E_2+Z$.
\end{enumerate}
\end{rem}

\begin{exa}[Bidouble covers of singular Del Pezzo surfaces of degree one]\label{exa:new}
See \cite{chennew}*{Section~2 and Section~3} for details.
Let $p_0, p_1, p_2, p_3$ be four points of $\PP^2$ in general position
and let $p_j'$ be the infinitely near point over $p_j$ corresponding to the line $\overline{p_0p_j}$ for $j=1, 2, 3$.
Finally, let $p$ be the eighth point satisfying the following Zariski open conditions:
\begin{enumerate}
\renewcommand{\labelenumi}{(\Roman{enumi})}
    \item $p \not \in \cup_{i=1}^3\{\overline{p_0p_i}:x_{i+1}=x_{i+2}\}
           \cup_{i=1}^3 \{ \overline{p_{i+1}p_{i+2}}: x_i=0\}$;
    \item $p \not \in c_1 \cup c_2 \cup c_3$, where $c_j$ is the unique conic passing through the five points
          $p_i, p_{i+1}, p_{i+1}', p_{i+2}$ and $p_{i+2}'$.
\end{enumerate}

\begin{figure}[ht]
\setlength{\unitlength}{1mm}
\centerline{\begin{picture}(30,25)
\linethickness{1pt}
\put(0,0){$p_1$}   \put(3,3){\circle*{1}}
\put(3,3){\vector(4,3){2}} \put(3,6){$p_1'$}
\put(3,3){\line(4,3){12}}
\put(28,0){$p_2$}  \put(27,3){\circle*{1}}
\put(27,3){\vector(-4,3){2}} \put(25,6){$p_2'$}
\put(27,3){\line(-4,3){12}}
\put(15,8){$p_0$}  \put(15,12){\circle*{1}}
\put(15,25){$p_3$} \put(15,22){\circle*{1}}
\put(15,22){\vector(0,-1){3}} \put(10,18){$p_3'$}
\put(15,22){\line(0,-1){10}}
\put(30,18){$p$}   \put(30,15){\circle*{1}}
\end{picture}}
\caption{\footnotesize{Configurations of $p_0, p_1,\ldots, p_3'$ and $p$} \label{figure2}}
\end{figure}
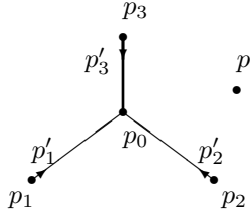


Denote by $E_j$ (respectively $E_j',$ $E$) the \textbf{total transform} of the point $p_j$ (respectively $p_j',$ $p$),
and by $L$ the pullback of a general line by $\sigma.$ Then
$\mathrm{Pic}(W)=\mathbb{Z}L\oplus \mathbb{Z}E_0\oplus \oplus_{j=1}^3 (\mathbb{Z}E_j\oplus \mathbb{Z}E_j')\oplus \mathbb{Z}E$
and $-K_W \equiv 3L-E_0-\sum\nolimits_{j=1}^{3}(E_j+E_j')-E.$
We list some properties of the surface $W$.
\begin{enumerate}[\upshape (1)]
    \item The surface $W$ is a weak Del Pezzo surface of degree one.
    \item There are exactly six nodal curves on $W.$ Their divisor classes are
          \begin{equation}\label{eq:(-2)curves}
          C_j \equiv L-E_0-E_j-E_j',\  C_j'\equiv E_j-E_j'\ \text{for}\ j=1, 2, 3.
          \end{equation}
          Let $\eta \colon W \rightarrow \Sigma$ be the morphism contracting there curves.
          The surface $\Sigma$ has six nodes and $-K_\Sigma$ is ample.
    \item The pencil of lines on $\PP^2$ passing
          through $p_0$ induces a fibration $g \colon W \rightarrow \PP^1$.
          Denote by $F$ a general fiber of $g$. Then $F \equiv L-E_0$.
          The fibration $g$ has exactly four singular fibers:
          $$\xymatrix @M=0pt@R=2pt{
          C_j              &   E_j'            &  C_j'   & \Gamma & E  \\
          \circ \ar@{-}[r] &  \circ \ar@{-}[r] &  \circ  & \circ \ar@{-}[r] &  \circ \\
          1&2&1&1&1}$$
          for $j=1, 2, 3$, where $\Gamma$ is the strict transform of the line $\overline{p_0p}$
          and $\Gamma \equiv L-E_0-E$.
    \item The linear system of $|-2K_W-\Gamma|$ consists of a $(-1)$-curve $\B_2$
          and the linear system of $|-2K_W-E|$ consists of a $(-1)$-curve $\B_3$.
          We have $\B_2\Gamma=\B_3E=3$ and $\B_2\B_3=\B_2E=\B_3\Gamma=1$.
          Moreover, the divisor $\Gamma+E+\B_2+\B_3$ has only nodes.
\end{enumerate}

Now we define three effective divisors on $W$
\begin{align}
    \Delta_1:&=F_b+\Gamma+(C_1+C_1'+C_2+C_2') \equiv 4L-4E_0-2E_1'-2E_2'-E,\nonumber\\
    \Delta_2:&=\B_2+(C_3+C_3') \equiv  -2K_W-2E_3'+E, \label{eq:newdata1}\\
    \Delta_3:&=\B_3 \equiv -2K_W-E. \nonumber
\end{align}
Here we require that the curve $F_b$ is a smooth fiber of $g$
and the divisor $\Delta:=\Delta_1+\Delta_2+\Delta_3$ has only nodes.
We also define three divisors
\begin{equation}\label{eq:newdata2}
    \begin{aligned}
     \L_1&=-2K_W-E_3',\\
     \L_2&=-K_W+(2L-2E_0-E_1'-E_2'-E),\\
     \L_3&=-K_W+(2L-2E_0-E_1'-E_2'-E_3').
    \end{aligned}
\end{equation}
It follows that $\Delta_i$ and $\L_i$ ($i=1, 2, 3$) satisfy \eqref{eq:buildingdata}.
These data define a smooth finite $G$-cover $\bpi \colon V \rightarrow W$
branched on the divisors $\Delta_1, \Delta_2$ and $\Delta_3$.
The (set theoretic) inverse image $\bpi^{-1}C_j$ or $\bpi^{-1}C_j'$ is a disjoint union of two $(-1)$-curves.
Let $\e \colon V  \rightarrow S$ be the blowdown of these twelve $(-1)$-curves.
From the construction, there is a finite $G$-cover $\pi \colon S \rightarrow \Sigma$ such that
the diagram \eqref{diag} commutes.
And $S$ is a smooth minimal surface of general type with $p_g(S)=0$ and $K_S^2=7$.
Moreover, $K_S$ is ample and the bicanonical map of $S$ is birational.
\end{exa}

In the notation of Section~2,
$\B_1=F_b+\Gamma, \N_1=C_1+C_1'+C_2+C_2', \N_2=C_3+C_3', \N_3=0$ and $D=-2K_W+\Gamma$.
Thus $(D\B_1, D\B_2, D\B_3)=(5, 5, 3)$ and $(\B_1\B_2, \B_1\B_3, \B_2\B_3)=(7, 5, 1)$.
So these surfaces satisfy \thmref{thm:numerical}~(b1) by \lemref{lem:DM}~(b).

We will need the following lemma to study deformations of surfaces in \exaref{exa:new}.
See Section~5.
\begin{lem}\label{lem:thethird-1curve}
Let $W$ be a weak Del Pezzo surface of degree one as in \exaref{exa:new}.
The linear system $|-K_W+E_3'-E|$ consists of a $(-1)$-curve $\Lambda$.
Moreover, $\Lambda$ intersects $F_b$ transversely and $\Lambda F_b=2$.
\end{lem}
\proof Note that $(-K_W+E_3'-E)^2=-1$ and $K_W(-K_W+E_3'-E)=-1$.
Also it is clear that $h^0(W, 2K_W-E_3'+E)=0$.
The Riemann-Roch theorem yields $h^0(W, \O_W(-K_W+E_3'-E)) \ge 1$.
By \eqref{eq:(-2)curves}, we can easily show that $-K_W+E_3'-E-C_j$ and $-K_W+E_3'-E-C_j'$ are not effective,
i.e., that any effective divisor $\Lambda$ in $|-K_W+E_3'-E|$ does not contain $C_j$ or $C_j'$ for $j=1, 2, 3$.
Because the nodal curves $C_1, \ldots, C_3'$ are exactly all the nodal curves of $W$
and $-K_W\Lambda=1$, we conclude that $\Lambda$ is irreducible and thus $\Lambda$ is a $(-1)$-curve.
We have $\Lambda F_b=(-K_W+E_3'-E)F_b=2$.
It remains to show that $\Lambda$ intersects $F_b$ transversely.

Note that $\Lambda C_k=\Lambda C_k'=0$ $(k=1, 2)$ and $\Lambda C_3=\Lambda C_3'=1$.
The linear system $|C_3+2\Lambda+C_3'|$ induces a genus $0$ fibration $g' \colon W \rightarrow \mathbb{P}^1$.
The same reasoning as the proof of \cite{chennew}*{Proposition~2.5}~(3) shows
$g'$ has exactly four singular fibers
          $$\xymatrix @M=0pt@R=2pt{
          (-2)             &   \Lambda_j           &  (-2)   & (-1) &  (-1) \\
          \circ \ar@{-}[r] &  \circ \ar@{-}[r] &  \circ  & \circ \ar@{-}[r] &  \circ \\
          1&2&1&1&1}$$
where the $(-2)$-curves are $C_1, \ldots, C_3'$ and $\Lambda_j$ is a $(-1)$-curve for $j=1, 2, 3$.
The curve $\Lambda$ is one of the curves $\Lambda_1, \Lambda_2$ and $\Lambda_3$.

Since $F_b(C_3+2\Lambda+C_3')=4$,
the restriction $g'|_{F_b} \colon F_b \rightarrow \mathbb{P}^1$ has degree $4$.
The curve $F_b$ is disjoint from the nodal curves, so $F_b\Lambda_j=2$ for $j=1, 2, 3$.
Because the multiplicity of $\Lambda_j$ in the singular fiber is $2$,
the ramification divisor of $g'|_{F_b}$ has degree at least $2 \times 3=6$,
and it has degree $6$ if and only if $F_b$ intersects $\Lambda_j$ transversely for $j=1, 2, 3$.
On the other hand, the Hurwitz formula implies that ramification divisor of $g'|_{F_b}$ has degree $6$.
We conclude that $F_b$ intersects $\Lambda_j$ transversely for $j=1, 2, 3$.
\qed

\section{Classification}
The whole section is devoted to prove \thmref{thm:classification}.
We have proved that $K_S$ is ample, $R_i^2=-1$ for $i=1, 2, 3$
and that $\Sigma$ is a rational surface containing $l_1+l_2+l_3$ nodes in \propref{prop:basicK2=7}.
So \thmref{thm:classification}~(c) follows from \thmref{thm:numerical}~(c).
We deal with one of the three cases (a), (b1) and (b2)  in \thmref{thm:numerical} in each subsection.
In Subsetion~4.1, we first show that the surfaces in the case (b1) are surfaces in \exaref{exa:new}.
This case is the most complicated one.
We give detailed proofs for this case and set up the classification process.
We then proceed analogously to exclude the case (b2) in Subseciton~4.2
and to show that the surfaces in case (a) are Inoue surfaces in Subsection~4.3.
This will complete the proof of \thmref{thm:classification}.

We make some general remarks and recall some general facts.
\begin{rem}\label{rem:recall}
\begin{enumerate}[\upshape (1)]
\item We consider the bidouble cover $\bpi \colon V \rightarrow W$ more often than the quotient map
      $\pi \colon S \rightarrow \Sigma$ (See the diagram \eqref{diag} for the notation).
      We remark that, on the smooth rational surface $W$,
      linear equivalence and numerical equivalence between divisors are the same.
\item We have $\B_i^2=-1$, $D\B_i=K_SR_i$ and $\B_i\B_{i+1}=R_iR_{i+1}$ for $i=1, 2, 3$. See \lemref{lem:DM}~(b).
\item The curve $\B_i$ is a disjoint union of irreducible smooth curves for $i=1, 2, 3$.
      The curve $\B_1+\B_2+\B_3$ is normal crossing and it is disjoint from the nodal curves
      $\N_1 \cup \N_2 \cup \N_3$. See \remref{rem:B}.
      The three curves $\B_1,\B_2,\B_3$ and
      the $l_1+l_2+l_3$ nodal curves $\N_1\cup \N_2 \cup \N_3$ generate a sublattice of $\mathrm{Num}(W)$.
      See \remref{rem:BBBNNN}.

\item The curve $\B_1+\B_2+\B_3$ does not contain any nodal curve;
      Otherwise, $W$ contains at least $l_1+l_2+l_3+1$ pairwise disjoint nodal curves.
      But this contradicts \propref{prop:basicK2=7}~(b), (d) and \cite{manynodes}*{Theorem~3.3}.
\item The divisor $M:=K_W+D$ is nef. See \lemref{lem:DM}~(g).
\end{enumerate}
\end{rem}

\subsection{Bidouble covers of Singular Del Pezzo Surfaces of Degree one}
In this subsection, we treat the case (b1) of \thmref{thm:numerical}.
Our aim is to prove that $S$ is a finite Galois $\mathbb{Z}_2^2$-cover of a singular Del Pezzo surface of degree one as described in \exaref{exa:new}.

In this case, we have
\begin{align}\label{eq:newDBBB}
    (D\B_1, D\B_2, D\B_3)=(5, 5, 3),\ (\B_1\B_2, \B_1\B_3, \B_2\B_3)=(7, 5, 1)\ \text{and}\ K_W^2=1
\end{align}
The minimal resolution $W$ of the quotient surface $\Sigma$ has
six disjoint nodal curves $\N_1 \cup \N_2$, where $\N_1$ consists of four nodal curves and
$\N_2$ consists of two nodal curves (Recall that $l_1=4, l_2=2$ and $l_3=0$).
We also have
\begin{align}
DK_W=\frac 12 D(D-\B_1-\B_2-\B_3)=-3,&& DM=D(K_W+D)=4, \label{eq:newDM}\\
K_WM=K_W(K_W+D)=-2,                  &&M^2=M(K_W+D)=2  \nonumber
\end{align}

The first lemma describes the surface $W$.
\begin{lem}\label{lem:553-KW}
The surface $W$ is a weak Del Pezzo surface of degree one.
\end{lem}
\proof
We have shown $K_W^2=1$. 
It suffices to show that $-K_W$ is nef.
Assume by contradiction that $-K_WC<0$ for an irreducible curve $C$.
The Riemann-Roch theorem implies $h^0(W, \O_W(-K_W)) \ge 2$.
So $C^2 <0$ and $C$ is contained in the fixed part of $|-K_W|$.
Because both $D$ and $M$ are nef, $1 \le  K_WC+DC=MC \le M(-K_W)=2$.
\lemref{lem:DM}~(f) implies $DC \ge 1$.
It follows that $K_WC=DC=1$ and $MC=2$.
So $C^2=-1$ or $-3$ by the adjunction formula.
On the other hand, because $M(-K_W-C)=0$ and $M^2=2$ by \eqref{eq:newDM},
the algebraic index theorem yields $(-K_W-C)^2=3+C^2 < 0$.
This gives a contradiction.
Hence $-K_W$ is nef and $W$ is a weak Del Pezzo surface of degree one.
\qed\\

The next lemma describes the branch divisors $\B_1, \B_2$ and $\B_3$.
\begin{lem}\label{lem:553BBB}The divisors $\B_1, \B_2, \B_3$ satisfy the following properties.
\begin{enumerate}[\upshape (a)]
\item The curve $\B_1$ is reducible: $\B_1=F_b+\Gamma$,
      where $\Gamma$ is a $(-1)$-curve and $\Gamma \equiv M+K_W$, while $F_b$ is a $0$-curve and $DF_b=4$.
\item The divisors $\B_2$ and $\B_3$ are $(-1)$-curves.
      Moreover, $\B_2 \equiv -2K_W-\Gamma$ and $\B_3 \equiv -2K_W-F_b+\Gamma$.
\end{enumerate}
\end{lem}
\proof
Because $\B_i^2=-1$ for $i=1, 2, 3$,
by \eqref{eq:DM} and  \eqref{eq:newDBBB}, we have
\begin{align*}
K_W\B_1=-3\ \text{and}\ K_W\B_2=K_W\B_3=-1.
\end{align*}
Note that
$-K_W$ has positive intersection numbers with each irreducible component of $\B_i$
by \remref{rem:recall}~(4).

Since $-K_W\B_1=3$ and $\B_1^2=-1$, the algebraic index theorem and the adjunction formula show that
$\B_1$  is a disjoint union of a $(-1)$-curve $\Gamma$ and a $0$-curve $F_b$.
Now we prove that $\Gamma \equiv M+K_W$.
Note that $\B_1(M+K_W)=\B_1(2K_W+D)=-1$.
It suffices to show that the linear system $|M+K_W|$ consists of a $(-1)$-curve.
Since $M$ is nef and big, it is $1$-connected.
The long exact sequence obtained from
$$0\rightarrow \O_W(K_W) \rightarrow \O_W(M+K_W) \rightarrow \omega_M \rightarrow0$$
gives $h^0(W, \O_W(M+K_W))=p_a(M)=1$.
Since $D(M+K_W)=1$, by \lemref{lem:DM}~(f), we may assume that $K_W+M \equiv \Phi+\Psi$,
where $\Phi$ is an irreducible curve with $D\Phi=1$ and $\mathrm{Supp}(\Psi) \in \N_1 \cup \N_2$.
So $\Phi^2 \le 0$ by the algebraic index theorem.
Since $K_W\Phi=K_W(M+K_W)=-1$, $\Phi$ is a $(-1)$-curve.
Then $\Psi^2=(M+K_W-\Phi)^2=(2K_W+D-\Phi)^2=0$ and thus $\Psi=0$.
Hence $|M+K_W|$ consists of a $(-1)$-curve and $\Gamma \equiv M+K_W$.
It follows that $D\Gamma=DM+DK_W=1$ and $DF_b=D(\B_1-\Gamma)=4$.
This proves (a).

For (b), because $-K_W\B_2=-K_W\B_3=1$, the curves $\B_2$ and $\B_3$ are irreducible.
Since $\B_2^2=\B_3^2=-1$, $\B_2$ and $\B_3$ are $(-1)$-curves.
Because $\Gamma \equiv M+K_W = 3K_W+\B_1+\B_2+\B_3$ and $\B_1 \equiv F_b+\Gamma$,
we have $\B_2+\B_3 \equiv -4K_W-F_b$.
So it suffices to show $\B_2 \equiv -2K_W-\Gamma$.

Note that $\B_2-(-2K_W-\Gamma) \equiv 3D-2\B_1-\B_2-2\B_3$.
It follows that $[\B_2-(-2K_W-\Gamma)]^2=0$ and $D[\B_2-(-2K_W-\Gamma)]=0$.
The algebraic index theorem yields $\B_2 \equiv -2K_W-\Gamma$.
\qed.\\

The next lemma describes a rational fibration on $W$.
See \cite{manynodes}*{Theorem~3.3}.
\begin{lem}\label{lem:553SF}
    The linear system $|F_b|$ induces a genus $0$ fibration $g \colon W \rightarrow  \mathbb{P}^1$.
    \begin{enumerate}[\upshape (a)]
        \item The curve $\Gamma$ is contained in a singular fiber $F_0$ of $g$ and $F_0=\Gamma+E$,
              where $E$ is a $(-1)$-curve and $E\Gamma=1$.
        \item The fibration $g$ has exactly four singular fibers:
              $F_0$ and $C_j+2E_j'+C_j'$ ($j=1, 2, 3$),
              where $E_j'$ is a $(-1)$-curve,
              while $C_j$ and $C_j'$ are the nodal curves contained in $\N_1 \cup \N_2$.
        \item If $E_0$ is a smooth section of $g$ and $E_0^2 <0$, then $E_0$ is a $(-1)$-curve.
    \end{enumerate}
\end{lem}
\proof
It is well known that a $0$-curve on a smooth rational surface induces a genus $0$ fibration.
Because $F_b\Gamma=0$, $\Gamma$ is contained in a fiber $F_0$ of $g$.
Similarly, the nodal curves $\N_1 \cup \N_2$ are contained in the fibers of $g$.
Denote by $F$ the general fiber of $g$.

Assume that $E$ is an irreducible component of $F_0$ such that $E\Gamma=1$.
Then $E$ is a smooth rational curve with $E^2<0$.
Since  $\B_2$ is a $(-1)$-curve and $\B_2F=(-2K_W-\Gamma)F=4$,
$\B_2 \not= E$ and thus $\B_2 E =(-2K_W- \Gamma) E \ge 0$.
It follows that $K_WE < 0$.
The adjunction formula shows that $E$ is $(-1)$-curve.
Then the Zariski lemma gives $F_0=\Gamma+E$.

Assume that $F_s$ is a singular fiber different from $F_0$ and $A$ is an irreducible component of $F_s$.
Then $A\Gamma=0$ since $A$ and $\Gamma$ are in different fibers.
Then $2K_WA=(\Gamma-D)A=-DA$.
By \propref{lem:DM}~(f),
$A$ is either a nodal curve in $\N_1 \cup \N_2$ or a $(-1)$-curve with $DA=2$.
Because $DF=4$, any singular fiber has one of the following types:
$$\xymatrix @M=0pt@R=2pt{
(-2)             &  (-1)          & (-2)       &(-1) & (-2)& (-1) & (-1) & (-1)\\
\circ \ar@{-}[r] &  \circ \ar@{-}[r] &  \circ  &\circ\ar@{-}[r] & \circ\ar@{-}[r] & \circ\ & \circ\ar@{-}[r] & \circ\\
1&2&1&1&1&1&1&1}$$
Note that each fiber of the first two types contributes $2$ to the Picard number $\rho(W)$.
The surface $W$ contains six disjoint nodal curves and $\rho(W)=10-K_W^2=9$.
We have seen that $g$ has one fiber $F_0=\Gamma+E$ of the third type.
So the other fibers are of the first type.

Assume that $E_0$ is a smooth section of $g$ with $E_0^2 <0$.
\lemref{lem:553-KW} implies $K_WE_0 \le 0$.
If $K_WE_0=0$, then $E_0$ is a nodal curve.
Also $E_0(\Gamma+E)=E_0F=1$.
By \lemref{lem:553BBB}~(b), either $\B_2E_0=-1$ or $\B_3E_0=-1$.
This is impossible because $\B_2$ and $\B_3$ are $(-1)$-curves.
Hence $K_WE_0<0$ and $E_0$ is a $(-1)$-curve by the adjunction formula.
\qed\\

\begin{lem}\label{lem:553nodalcurves}
The surface $W$ contains exactly six nodal curves $\N_1 \cup \N_2$.
\end{lem}
\proof
Assume that $C$ is a nodal curve different from the six nodal curves $\N_1 \cup \N_2$.
According to \lemref{lem:553SF}~(b), $C$ is not contained in the fibers of $g$.
\lemref{lem:553SF}~(c) implies  $FC \ge 2$.
Then, by \lemref{lem:553BBB}~(b), either $\B_2C<0$ or $\B_3C<0$.
This contradicts that $\B_2$ and $\B_3$ are $(-1)$-curves.
Thus $W$ contains exactly six nodal curves $\N_1 \cup \N_2$.
\qed\\

We have shown in diagram \eqref{diag} that the cover $\bpi \colon V \rightarrow W$ is branched on
$\B_1+\N_1,  \B_2+\N_2$ and $\B_3$, where
\begin{align}
\B_1 =F_b+\Gamma,&& \B_2\equiv -2K_W-\Gamma,&& \B_3 \equiv -2K_W-E \equiv -2K_W-F+\Gamma \label{eq:553Bclass}
\end{align}
Now we consider how the nodal curves $C_1, \ldots, C_3'$ in the fibers of $g$ (see \lemref{lem:553SF})
distribute along two divisors $\N_1$ and $\N_2$.
We conclude that $\N_1$ and $\N_2+F$ are divisible by $2$ in $\mathrm{Pic}(W)$
from \eqref{eq:buildingdata} and \eqref{eq:553Bclass}.
Since $C_i+C_i' \equiv F-2E_i'$ for $i=1, 2, 3$,
we may assume
\begin{align}
    \N_1=(C_1+C_1')+(C_2+C_2')\equiv 2F-2E_1'-2E_2',\ \N_2=C_3+C_3'\equiv F-2E_3'. \label{eq:553Nclass}
\end{align}

Finally, we shall show that $W$ arises as the blowup of $\mathbb{P}^2$ as described in \exaref{exa:new}.
We first claim that there exists a smooth section $E_0$ of $g$,
such that $E_0$ is a $(-1)$-curve, $E_0\Gamma=1$ and $E_0E=0$.
Let $p \colon W \rightarrow W'$ be the morphism blowing down $E$.
Because $K_{W'}^2=2$, the fiberation $g' \colon W \rightarrow \mathbb{P}^1$ induced by $g$ is not relatively minimal.
Thus $g'$ has a smooth section $E_0'$ such that $E_0'^2<0$.
The strict transform $E_0$ of $E_0'$ is a smooth section of $g$ and thus it is a $(-1)$-curve by \lemref{lem:553SF}~(c).
It follows that $E_0\Gamma=1$ and $E_0 E=0$.
The claim is proved.

Since $F \equiv C_j+2E_j'+C_j'$, we have $E_0E_j'=0$.
After possibly relabeling $C_j$ and $C_j'$ for each $j=1, 2, 3$,
we may assume that $E_0C_j=1$ and $E_0C_j'=0$ .

Blowing down $E_j' ( j=1, 2, 3)$, then blowing down the image of $C_j' ( j=1, 2, 3)$ and finally blowing down
$E_0$ and $E$, we obtain a birational morphism $\sigma \colon W \rightarrow \mathbb{P}^2$.
Let $p_0:=\sigma(E_0), p_1:=(E_j' \cup C_j')\ (j=1, 2, 3)$ and $p :=\sigma(E)$.
Denote by $p_j'$ the infinitely near point over $p_j$ corresponding to the line $\overline{p_jp_0}$.
Note that the exceptional divisor on $W$ corresponding to $p_j'$ is $E_j'$.
The fibration $g$ corresponds to the pencil of lines passing through the point $p_0$.
Then \lemref{lem:553SF}~(b) implies that $p_0, p_i, p_{i+1}$ are not collinear
and $p_0, p, p_i$ are not collinear for each $i=1, 2, 3$.
\lemref{lem:553nodalcurves} implies that $p_1, p_2, p_3$ are not collinear and the point $p$
satisfies  the conditions (I) and (II) in \exaref{exa:new}.
Otherwise, for example, if $p \in c_1$, then the strict transform of $c_1$ on $W$ is a nodal curve,
which is different from $C_1, \ldots, C_3'$, a contradiction to \lemref{lem:553nodalcurves}.

Hence we have shown that $S$ is a surface in \exaref{exa:new}.

\subsection{Exclusion}
Our next goal is to exclude the case (b2) of \thmref{thm:numerical}.
We assume by contradiction that there is a pair $(S, G)$ satisfying the property (b2) of \thmref{thm:numerical}.
In this case, we have
\begin{align}\label{eq:exclusionDBBB}
(D\B_1, D\B_2, D\B_3)=(5, 5, 3),\ (\B_1\B_2, \B_1\B_3, \B_2\B_3) \in \{(3, 5, 1), (7, 1, 1)\}\ \text{and}\ K_W^2=-1
\end{align}
The minimal resolution $W$ of the quotient surface $\Sigma$ has eight pairwise disjoint nodal curves
$\N_1 \cup \N_2 \cup \N_3$.
We also have
\begin{align}
DK_W=\frac 12 D(D-\B_1-\B_2-\B_3)=-3, && D M=D(K_W+D)=4, \label{eq:exclusionDM}\\
K_WM=K_W(K_W+D)=-4,                   && M^2=M(K_W+D)=0. \nonumber
\end{align}

\begin{lem}\label{lem:exclusionrationalfibration}
\begin{enumerate}[\upshape (a)]
\item The linear system $|M|$ is composed with a pencil $|F|$ and $M \equiv 2F$,
      where $|F|$ is a base point free pencil of rational curves and $DF=2$.
\item Let $g \colon W \rightarrow \PP^1$ be the fibration defined by $|F|$.
      Then $g$ has exactly five singular fibers: $E_1+E_2$ and $C_j+2\Gamma_j+C_j'$ $(j=1,2,3,4)$,
      where $E_1$, $E_2$ and $\Gamma_j$ are $(-1)$-curves,
      while $C_j$ and $C_j'$ are the nodal curves in $\N_1 \cup \N_2 \cup \N_3$.
      Moreover, $DE_1=DE_2=1$.
\end{enumerate}
\end{lem}
\proof
\lemref{lem:DM}~(e) gives $h^0(W, \O_W(M))=p_a(D)=3$.
Assume $|M|=|\Phi|+\Psi$, where $|\Phi|$ is the moving part and $\Psi$ is the fixed part.
Because both $M$ and $\Phi$ is nef, and $M^2=0$,  we have $M\Phi=M\Psi=0$ and $\Phi^2=\Phi\Psi=\Psi^2=0$.
Then $|\Phi|$ is composed a pencil $|F|$ and $\Phi \equiv 2F$.
\lemref{lem:xiao} implies $DF \ge 2$.
Since $DM=4$ and $D\Phi=4$, we have $D\Psi=0$.
Then \lemref{lem:DM}~(f) and $\Psi^2=0$ imply $\Psi=0$.
Thus $M \equiv \Phi \equiv 2F$, $DF=2$ and $K_WF=-2$.
Therefore $|F|$ is a base point free pencil of rational curves.

Since $F\N_i=\frac 12 (D+K_W)\N_i=0$ for $i=1, 2, 3$,
the nodal curves $\N_1 \cup \N_2 \cup \N_3$ are contained in the singular fibers of $g$.
Assume that $A$ is an irreducible component of a singular fiber of $g$.
Then $A$ is a smooth rational curve with $A^2<0$.
Note that $0=2FA=MA=K_WA+DA \ge K_WA=-2-A^2$.
By \lemref{lem:DM}~(f),  $A$ is either one of the nodal curves in $\N_1 \cup \N_2 \cup \N_3$
or a $(-1)$-curve with $DA=1$.
The rest of the proof is similar to that of \lemref{lem:553SF}~(a)-(b).
\qed

\begin{lem}\label{lem:exclusionBBB}
\begin{enumerate}[\upshape (a)]
\item The divisors $\B_1, \B_2$ and $\B_3$ satisfy the following properties:
      $\B_1\B_2=7, \B_1\B_3=\B_2\B_3=1$ and $\B_iF=2$ for $i=1, 2, 3$.
\item The divisor $\B_3$ is an irreducible smooth elliptic curve and $\B_3 \equiv -K_W$.
\end{enumerate}
\end{lem}
\proof
Note that $F \nequiv =\frac 12(D+K_W) \nequiv=\frac 14(3D-\B_1-\B_2-\B_3)$ and $F\N_i=0$ for $i=1, 2, 3$.
By \eqref{eq:buildingdata}, $F\B_1, F\B_2$ and $F\B_3$ are of the same parity.
If $(\B_1\B_2, \B_1\B_3, \B_2\B_3)=(3, 5, 1)$,
then $F\B_1=2$ and $F\B_2=3$.
This gives a contradiction.
Hence $\B_1\B_2=7$ and $\B_1\B_3=\B_2\B_3=1$.
It follows that $F\B_i=2$ for $i=1, 2, 3$.

By \remref{rem:BBBNNN},  from the intersection numbers $D\B_i$  and $D\N_i=0$,
we can write $D$ as a $\mathbb{Q}$-linear combination of $\B_1, \B_2$ and $\B_3$.
It turns out to be a $\mathbb{Z}$-linear combination: $D \equiv \B_1+\B_2-\B_3$.
Since $\mathrm{Pic}(W)$ has no torsion, $\B_3 \equiv -K_W$ by \eqref{eq:DM}.
Assume that $\B_3 =\Phi_1+\ldots+\Phi_t$, where the curves $\Phi_k$ are irreducible, smooth and pairwise disjoint.
Then $\Phi_1^2=\Phi_1\B_3=-K_W\Phi_1$. It follows that $\Phi_1$ is a smooth elliptic curve.
The long exact sequence obtained from
$$0 \rightarrow \O_W(K_W) \rightarrow \O_W(K_W+\Phi_1) \rightarrow \omega_{\Phi_1} \rightarrow 0$$
gives $h^0(W, \O_W(K_W+\Phi_1))=1$.
Since $K_W+\Phi_1 \equiv -\Phi_2-\ldots-\Phi_t$, $t=1$ and $\B_3=\Phi_1$.
This is the desired conclusion.
\qed

\begin{lem}\label{lem:exclusionellitpicfibration}
The linear systems $|\B_3+E_1|$ and $|\B_3+E_2|$ are base point free pencil of elliptic curves.
\end{lem}
\proof

Note that $(\B_1+\B_2)E_k=(D-K_W)E_k=(M-2K_W)E_k=2$ for $k=1, 2$.
Since $E_k$ is disjoint from the nodal curves $\N_1 \cup \N_2 \cup \N_3$, by \eqref{eq:buildingdata},
$\B_1E_k, \B_2E_k$ and $\B_3E_k$ are of the same parity.
Since $\B_3E_k=-K_WE_k=1$, $\B_1E_k$ and $\B_2E_k$ are odd integers.

Fix $k \in \{1, 2\}$.
Assume $\B_1+\B_2 \not \ge E_k$.
Then $\B_1E_k=\B_2E_k=1$.
The intersection number matrix of $\B_1, \B_2, \B_3$ and $E_k$ is nondegenerate.
This contradicts \remref{rem:BBBNNN}.
Therefore $\B_1+\B_2 \ge E_k$ for $k=1, 2$.

Without loss of generality, we assume $\B_1 \ge E_1$.
Then $\B_1E_1=-1$ and $\B_2E_1=3$.
By \remref{rem:BBBNNN},
we may write $E_1$ as the $\mathbb{Q}$-linear combination of $\B_1, \B_2, \B_3$.
That is $E_1 \nequiv \frac 12 (\B_1-\B_3)$, i.e., $\B_1 \equiv 2E_1+\B_3$.
Then $\B_1E_2=(2E_1+\B_3)E_2=3$ and $\B_2E_2=2-\B_1E_2=-1$.
Thus $\B_2 \ge E_2$ and $\B_2 \equiv 2E_2+\B_3$ again by \remref{rem:BBBNNN}.

For $k=1, 2$,
the divisor $\B_k$ is a disjoint union of smooth irreducible curves, so is $\B_k-E_k$.
Also $\B_k-E_k$ and $\B_3+E_k$ have no common irreducible component (cf.~\remref{rem:B}).
Since $\B_k \equiv \B_3+2E_k$, $\B_k-E_k$ and $\B_3+E_k$ generate a base point free pencil of elliptic curves
for $k=1, 2$.\qed\\

We are now in a position to deduce a contradiction.
Let $P:=E_1+E_2+\B_3$.
We have proved that $E_1E_2=\B_3E_1=\B_3E_2=1$ and $\B_3, E_1, E_2$ have no common points (cf.~\remref{rem:B}).
Because $|E_1+E_2|, |\B_3+E_1|$ and $|\B_3+E_2|$ are base point free pencils of curves,
$|P|$ is base point free.
Let $F_1$ be a general element of $|\B_3+E_1|$.
Then there is an exact sequence
$$0 \rightarrow \O_W(E_2) \rightarrow \O_W(P) \rightarrow \O_{F_1}(P) \rightarrow 0.$$
It is easy to see that $h^0(W, \O_W(E_2))=1$ and $h^1(W, \O_W(E_2))=0$.
Since $F_1$ is smooth elliptic curve and $PF_1=2$, we have $h^0(W, \O_{F_1}(P))=2$.
Then $h^0(W, \O_W(P))=3$ and the trace of $|P|$ on $F_1$ is complete.
The linear system $|P|$ defines a morphism $f \colon W \rightarrow \PP^2$ of degree $P^2=3$.
However, the restriction of $f$ on a general element $F_1$ in $|\B_3+E_1|$ has degree $PF_1=2$,
which was supposed to be a factor of $\deg f=3$.
This gives a contradiction and we exclude the case (b2) in \thmref{thm:numerical}.

\subsection{Inoue Surfaces}
The subsection is devoted to classify the pairs $(S, G)$ satisfying the property (a) of \thmref{thm:numerical}.
We show that $S$ is an Inoue surface in \exaref{exa:inoue}.

In this case, we have
\begin{align}\label{eq:inoueDBBB}
(D\B_1, D\B_2, D\B_3)=(7, 5, 5), (\B_1\B_2, \B_1\B_3, \B_2\B_3)=(5, 9, 7)\ \text{and}\ K_W^2=3
\end{align}
The minimal resolution $W$ of the quotient surface $\Sigma$ has
four disjoint nodal curves $\N_1 \cup \N_3$, where $\N_1$ consists of two nodal curves and so does $\N_3$
(Recall that $l_1=2, l_2=0$ and $l_3=2$).
We also have
\begin{align}
DK_W=\frac 12 D(D-\B_1-\B_2-\B_3)=-5, && DM=D(K_W+D)=2,\label{eq:inoueDM}\\
K_WM=K_W(K_W+D)=-2,                   && M^2=M(K_W+D)=0. \nonumber
\end{align}

The following lemma describes the surface $W$.
\begin{lem}\label{lem:755-KW}
    \begin{enumerate}[\upshape (a)]
        \item The linear system $|M|$ is a base point free pencil of rational curves.
        \item The surface $W$ is a weak Del Pezzo surface of degree three.
    \end{enumerate}
\end{lem}
\proof
The linear system  $|M|$ is a pencil of curves by \propref{lem:DM}~(e).
Since $M^2=0$ by \eqref{eq:inoueDM}, it suffices to prove that $M$ is base point free.
The proof is similar to that of \lemref{lem:exclusionrationalfibration}~(a).

The Riemann-Roch Theorem implies $h^0(W, \O_W(-K_W)) \ge 4$.
To prove (b), it suffices to prove that $-K_W$ is nef.
The rest of the proof runs as that of \lemref{lem:553-KW}.
\qed\\

The next lemma determines the branch divisors $\B_1, \B_2$ and $\B_3$.
\begin{lem}\label{lem:755BBB}The divisors $\B_1, \B_2$ and $\B_3$ are as follows:
\begin{align}\label{eq:755branchdivisor}
\B_1=\Gamma_1+F_2,&& \B_2=\Gamma_2+F_3,&& \B_3=\Gamma_3+F_1+F_1',
\end{align}
where the curves $\Gamma_i$ are $(-1)$-curves, $F_i$ and $F_1'$ are $0$-curves,
$\Gamma_i+F_i \equiv -K_W$ for $i=1, 2, 3$
and $\Gamma_1+\Gamma_2+\Gamma_3 \equiv -K_W$.
Moreover, $F_1, F_1' \in |M|$.
\end{lem}
\proof
By \eqref{eq:DM} and \eqref{eq:inoueDBBB}, we have
\begin{align*}
K_W\B_1=K_W\B_2=-3\ \text{and}\ K_W\B_3=-5.
\end{align*}
Because $M\B_3=(D+K_W)\B_3=0$,
$\B_3$ is contained in some divisors of the linear system $|M|$.
The Zariski's lemma, the algebraic index theorem and the adjunction formula imply
$\B_3=\Gamma_3+F_1+F_1'$, where $\Gamma_3$ is $(-1)$-curve,
while $F_1$ and $F_1'$ are two $0$-curves in $|M|$ (See \remref{rem:recall}~(3) and (4)).
The same reasoning applies to $\B_2$ and $\B_3$ yields
$\B_1=\Gamma_1+F_2$ and $\B_2=\Gamma_2+F_3$,
where $\Gamma_1$ and $\Gamma_2$ are $(-1)$-curves, while $F_2$ and $F_3$ are $0$-curves.

We claim that $\Gamma_i\Gamma_{i+1}=1$ for $i=1, 2, 3$.
Actually, since $-K_W(\Gamma_i+\Gamma_{i+1})=2$,  the algebraic index theorem implies $\Gamma_i\Gamma_{i+1} \le 1$.
As an irreducible component of the curve $\B_i$,
$\Gamma_i$ is disjoint from the nodal curves $\N_1 \cup \N_3$.
Now fix $i$.
If $\Gamma_i\Gamma_{i+1}=0$,
then blowing down $\Gamma_i$ and $\Gamma_{i+1}$,
we obtain a rational surface $W'$ containing four disjoint nodal curves and $\rho(W')=5$.
This gives a contradiction to \cite{manynodes}*{Theorem~3.3}.
The claim is proved and thus $(\Gamma_1+\Gamma_2+\Gamma_3)^2=3$.
Since $-K_W(\Gamma_1+\Gamma_2+\Gamma_3)=3$, the algebraic index theorem implies $\Gamma_1+\Gamma_2+\Gamma_3 \equiv-K_W$.

Since $F_1, F_1' \in |M|$ and $M \equiv K_W+D \equiv 3K_W+\B_1+\B_2+\B_3$,
$F_1+F_2+F_3 \equiv -2K_W$ by \eqref{eq:755branchdivisor}.
It follows that $F_i(F_1+F_2+F_3)=4$ and $F_iF_{i+1}=2$ for $i=1, 2, 3$.

Since $\B_i$ is a disjoint union of smooth curves, we have $\Gamma_1F_2=0, \Gamma_2F_3=0$ and $\Gamma_3F_1=0$.
Also $5=\B_1\B_2=(\Gamma_1+F_2)(\Gamma_2+F_3)=3+\Gamma_1F_3+\Gamma_2F_2$, i.e., $\Gamma_1F_3+\Gamma_2F_2=2$.
Similarly, $\B_1\B_3=9$ and $\B_2\B_3=7$ imply that $2\Gamma_1F_1+\Gamma_3F_2=4$ and $2\Gamma_2F_1+\Gamma_3F_3=2$.
Note that $(F_1+F_2+F_3)\Gamma_i=-2K_W\Gamma_i=2$ for $i=1, 2, 3$.
It follows that $\Gamma_iF_i=2$ and $F_i \Gamma_{i+1}=F_i\Gamma_{i+2}=0$ for $i=1, 2, 3$.
Since $(\Gamma_i+F_i)^2=3$ and $-K_W(\Gamma_i+F_i)=3$, the algebraic index theorem gives $\Gamma_i+F_i\equiv -K_W$.
\qed\\

\begin{lem}\label{lem:755nodalcurves}
The surface $W$ is the minimal resolution of the $4$-nodal cubic surface.
\end{lem}
\proof
It is well known that any weak Del Pezzo surface of degree three
is the minimal resolution of a normal cubic surface in $\mathbb{P}^3$.
We have seen that $W$ contains four nodal curves $\N_1 \cup \N_3$.
Assume that $C$ is a nodal curve on $W$.
Because $\Gamma_1+\Gamma_2+\Gamma_3 \equiv -K_W$ and $\Gamma_i+F_i \equiv -2K_W$,
$\Gamma_iC=F_iC=0$ for $i=1, 2, 3$.
Then $DC=0$ by \eqref{eq:DM} and \eqref{eq:755branchdivisor}.
\lemref{lem:DM}~(f) implies $C \in \N_1 \cup \N_3$.
Hence $W$ contains exactly four nodal curves and it is the minimal resolution of the $4$-nodal cubic surface.
\qed\\

We have shown in diagram \eqref{diag} that the cover $\bpi \colon V \rightarrow W$ is branched on $\B_1+\N_1,  \B_2$
and $\B_3+\N_3$, where both $\N_1$ and $\N_3$ consist of two disjoint nodal curves,
$\B_1, \B_2$ and $\B_3$ are described in \lemref{lem:755BBB}.
We conclude that $F_2+\N_1$ and $F_2+\N_3$ are divisible in $\mathrm{Pic}(W)$ by
\eqref{eq:buildingdata} and \eqref{eq:755branchdivisor}.
It follows that the two nodal curves in $\N_k$ ($k=1, 3$) are in the same singular member of the pencil $|F_2|$.
Comparing with \eqref{eq:755branchdivisor} and \exaref{exa:inoue},
we conclude that $S$ is an Inoue surface by \remref{rem:inouerem}.

\section{Deformations and the moduli space}
We study the local deformations and the moduli of the surfaces in the case (b) of \thmref{thm:classification}.
These surfaces are exactly the ones in \exaref{exa:new}.
Our goal is to prove \thmref{thm:moduli}.
Throughout this section, we assume that $S$ is a smooth minimal surface in \exaref{exa:new}.
We denote by $\t_S$ the tangent sheaf of $S$.
The following proposition estimates the dimension of the cohomology group $H^2(S, \t_S)$.

\begin{prop}\label{prop:eigenspaces}
The dimensions of the eigenspaces of $H^2(S, \t_S)$ ( for the $G$-action) satisfy the following properties:
$$\dim H^2(S, \t_S)^{\mathrm{inv}}=0,\ \dim H^2(S, \t_S)^{\mathrm{\chi_1}} \le 2,\ \dim H^2(S, \t_S)^{\mathrm{\chi_2}}\le 2,\ \dim H^2(S, \t_S)^{\mathrm{\chi_3}} \le 3.$$
\end{prop}
We use the methods in \cite{burniatII},  \cite{inouemfd}, \cite{burniatIII}
and \cite{chenKNMP} to prove \propref{prop:eigenspaces}.
The techniques involved depend on the exact sequences in \cite{ev}*{Properties~2.3~(c), p.~13}.
We also need the following lemma, which generalizes \cite{burniatII}*{Lemma~5.1}.
\begin{lem}[\cite{chenKNMP}*{Lemma~4.4}]\label{lem:logsheaf}
Let $X$ be a projective smooth surface.
Let $Y_1, \ldots, Y_{k-1}$ and $Y$ be $k$ irreducible smooth curves on $X$.
Assume that $Y_1+\ldots+Y_{k-1}+Y$ has only nodes.
Then there is an exact sequence
\begin{align*}
0 &\rightarrow \o_{X}^1(\log Y_1, \ldots, \log Y_{k-1}, \log Y) \rightarrow
\o_{X}^1(\log Y_1, \ldots, \log Y_{k-1})(Y) \\
&\rightarrow \o_{Y}^1(Y_1+\ldots+Y_{k-1}+Y) \rightarrow 0.
\end{align*}
Moreover, if $\L$ is a divisor of $X$ and
$Y.(K_{X}+2Y+Y_1+\ldots+Y_{k-1}+\L) <0,$ then
\begin{align*}
\dim H^0(X, \o_{X}^1(\log Y_1, \ldots, \log Y_{k-1})(Y+\L))=
\dim H^0(X, \o_{X}^1(\log Y_1, \ldots, \log Y_{k-1}, \log Y)(\L)).
\end{align*}
\end{lem}

\proof[Proof of \propref{prop:eigenspaces}]
We recall in the \exaref{exa:new} that $V$ is a blowup of the surface $S$
and $\bpi \colon V \rightarrow W$ is a finite Galois $G$ $(\cong \mathbb{Z}_2^2)$-cover branched on
$\Delta_1, \Delta_2$ and $\Delta_3$ (see the diagram \eqref{diag} and \eqref{eq:newdata1}).
We often refer to \exaref{exa:new} (2)-(4) for the intersection numbers of the curves and
the classes of curves in the Picard group $\mathrm{Pic}(W) \cong H^2(W, \mathbb{Z})$.

Denote by $\t_V$ the tangent sheaf of the surface $V$.
Because blowing down a $(-1)$-curve does not change the dimension of the second cohomology
group of the tangent sheaf,
we have $\dim H^2(S,\t_S)^{\mathrm{inv}}= \dim H^2(V,\t_{V})^{\mathrm{inv}}$ and
$\dim H^2(S,\t_V)^{\chi_i}=\dim H^2(V,\t_{V})^{\chi_i}$
for $i=1,2,3$.
Serre duality implies $H^2(V, \t_{V})\cong H^0(V, \o^1_{V}\otimes\o^2_{V})^*$.
The $G$-cover structure gives (cf.~\cite{bidoublecover}*{Theorem~2.16})
\begin{align}
\bpi_{\ast}(\o^1_{V} \otimes \o^2_{V})&=\o^1_{W}(\log \Delta_1, \log \Delta_2, \log \Delta_3)(K_W)
\oplus (\oplus_{i=1}^3 \o^1_{W}(\log \Delta_i)(K_{W}+\L_i)) \label{eq:eignesheaves}\\
H^2(V,\o^1_{V}\otimes\o^2_{V})&= H^0(W,\o^1_{W}(\log \Delta_1, \log \Delta_2, \log \Delta_3)(K_W))\nonumber\\
                              &\oplus (\oplus_{i=1}^3 H^0(W,\o^1_{W}(\log \Delta_i)(K_{W}+\L_i))), \nonumber
\end{align}
where the invertible sheaves $\L_1, \L_2, \L_3$ are given by \eqref{eq:newdata2}.
It is sufficient to calculate the dimension of each summand.\\

We first prove that $H^0(W,\o^1_{W}(\log \Delta_1, \log \Delta_2, \log \Delta_3)(K_W))=0$.
We have the following exact sequence from \cite{bidoublecover}*{Lemma~3.7} and \cite{maximum}*{Lemma~3}
\begin{equation}\label{eq:newinv.log.seq}
0 \rightarrow \o_{W}^1(K_{W}) \rightarrow
          \o_{W}^1(\log \Delta_1, \log \Delta_2, \log \Delta_3)(K_{W}) \rightarrow
          \oplus_{i=1}^3\O_{\Delta_i}(K_{W}) \rightarrow
          0
\end{equation}

Since $H^0(W, \o_{W}^1)=0$ and $-K_{W}$ is effective,
$H^0(W, \o_{W}^1(K_{W}))=0$.
It suffices to show that the boundary map
$$\delta \colon H^0(W,\oplus_{i=1}^3 \O_{\Delta_i}(K_{W}))\rightarrow H^1(W,\o_{W}^1(K_{W}))$$
is injective.
The linear system $|-K_W|$ has only one simple base point because $W$ is a weak Del Pezzo surface of degree one.
So there is a morphism $\O_{W}(K_{W}) \rightarrow \O_{W}$,
which is not identically zero on any component of the divisors $\Delta_i$.
This morphism $\O_{W}(K_{W}) \rightarrow \O_{W}$ induces a commutative diagram
$$\xymatrix@C=1.5em@R=1.5em{
  0\ar[r] & \o_{W}^1(K_{W}) \ar[d] \ar[r]
          & \o_{W}^1(\log \Delta_1, \log \Delta_2, \log \Delta_3)(K_{W}) \ar[d] \ar[r]
          & \oplus_{i=1}^3\O_{\Delta_i}(K_{W}) \ar[d] \ar[r]
          & 0\\
  0\ar[r] & \o_{W}^1 \ar[r]
          & \o_{W}^1(\log \Delta_1, \log \Delta_2, \log \Delta_3) \ar[r]
          & \oplus_{i=1}^3\O_{\Delta_i} \ar[r]
          & 0
}$$
It gives a commutative diagram of cohomology groups
$$\xymatrix@R=1.5em{
  H^0(W,\oplus_{i=1}^3 \O_{\Delta_i}(K_{W}))\ \ar[d]_{\psi_2} \ar[r]^-{\delta} \ar[1,1]^{\psi}
& H^1(W, \o_{W}^1(K_{W})) \ar[d]\\
  H^0(W, \oplus_{i=1}^3 \O_{\Delta_i}) \ar[r]^-{\psi_1}
& H^1(W, \o_{W}^1)
}$$

Note that
$H^0(W,\oplus_{i=1}^3 \O_{\Delta_i}(K_{W})) \cong \oplus_{i=1}^3 H^0(W, \O_{C_i}) \oplus H^0(W, \O_{C_i'})$.
By \cite{bidoublecover}*{Lemma~3.7},
the image of the function identically equal to $1$
on $C_i$ (respectively $C_i'$) maps under $\psi_1$ to
the first Chern class of $C_i$ (respectively $C_i'$).
Because the curves $C_i$ and $C_i'$ are disjoint nodal curves,
their Chern classes are linearly independent in $H^1(W,\o^1_{W})$.
Thus the composite map $\psi=\psi_1\psi_2$ is injective.
It follows that $\delta$ is also injective.
We thus get
$$\dim H^2(S, \t_S)^{\mathrm{inv}}=\dim H^0(W, \o^1_{W}(\log \Delta_1, \log \Delta_2, \log \Delta_3)(K_W))=0.$$

We  now calculate $\dim H^0(W, \o_W^1(\log \Delta_1)(K_W+\L_1))$.
Applying the last statement of \lemref{lem:logsheaf} to $\Gamma, C_3$ and $C_3'$, we have
\begin{align}
\lefteqn{\dim H^0(W, \o_W^1(\log \Delta_1)(K_W+\L_1))=} \label{eq:chi11}\\
& &\dim H^0(W, \o_W^1(\log F_b, \log C_1, \log C_1', \log C_2, \log C_2', \log C_3, \log C_3')(K_W+\L_1-C_3-C_3'+\Gamma)) \nonumber
\end{align}
Note that $K_W+\L_1-C_3-C_3'+\Gamma \equiv (-K_W-E_3')-(L-E_0-2E_3')+(L-E_0-E)=-K_W+E_3'-E$
(see \eqref{eq:newdata2} and \eqref{eq:(-2)curves}).
By \lemref{lem:thethird-1curve},
$|-K_W+E_3'-E|$ consists of a $(-1)$-curve $\Lambda$ and $\Lambda$ meets $F_b$ transversely.
Now $\Lambda(K_W+2\Lambda+F_b+C_1+\ldots+C_3')=1$,
\lemref{lem:logsheaf} gives the following exact sequence
\begin{align*}
\lefteqn{0 \rightarrow \o_W^1(\log F_b, \log C_1, \log C_1', \log C_2, \log C_2', \log C_3, \log C_3', \log \Lambda)} \\
& &\rightarrow \o_W^1(\log F_b, \log C_1, \log C_1', \log C_2, \log C_2', \log C_3, \log C_3')(\Lambda)
\rightarrow \O_\Lambda(1) \rightarrow 0.
\end{align*}
The first Chern classes of $F_b, C_1, \ldots, C_3'$ and $\Lambda$ are linearly independent in $H^2(W, \mathbb{C})$.
By \cite{bidoublecover}*{Lemma~3.7}, it follows that
$$H^0(W, \o_W^1(\log F_b, \log C_1, \log C_1', \log C_2, \log C_2', \log C_3, \log C_3', \log \Lambda))=0.$$
From the long exact sequence of cohomology we conclude that
\begin{align}
\dim H^0(W, \o_W^1(\log F_b, \log C_1, \log C_1', \log C_2, \log C_2', \log C_3, \log C_3')(\Lambda)) \le 2 \label{eq:chi12}
\end{align}
We thus get
$\dim H^2(S, \t_S)^{\mathrm{\chi_1}}=h^0(W, \o_W^1(\log \Delta_1)(K_W+\L_1)) \le 2$ by \eqref{eq:chi11} and \eqref{eq:chi12}.\\

We proceed to calculate $\dim H^0(W, \o_W^1(\log \Delta_2)(K_W+\L_2))$.
Note that $\Gamma$ intersects $\B_2$ transversely (see \exaref{exa:new}~(4))
and $\Gamma (K_W+2\Gamma+\B_2+C_3+C_3'+(K_W+\L_2-\Gamma))=0$.
\lemref{lem:logsheaf} yields an exact sequence
\begin{align*}
\lefteqn{0 \rightarrow \o_W^1(\log \B_2, \log C_3, \log C_3',  \log \Gamma)(K_W+\L_2-\Gamma)}\\
  & & \rightarrow \o_W^1(\log \B_2, \log C_3, \log C_3')(K_W+\L_2) \rightarrow \O_{\Gamma} \rightarrow 0.
\end{align*}
It follows that
\begin{align}
\lefteqn{\dim H^0(W, \o_W^1(\log \B_2, \log C_3, \log C_3')(K_W+\L_2))} \nonumber \\
& &\le \dim H^0(W, \o_W^1(\log \B_2, \log C_3, \log C_3', \log \Gamma)(K_W+\L_2-\Gamma))+1 \label{eq:chi21}
\end{align}
Note that $\B_2(K_W+\L_2-\Gamma)=0$.
Tensoring the following exact sequence (cf.~\cite{ev}~Properties~2.3~(b))
\begin{align*}
\lefteqn{0 \rightarrow \o_W^1(\log C_3, \log C_3', \log \Gamma)} \nonumber\\
& &\rightarrow \o_W^1(\log \B_2, \log C_3, \log C_3', \log \Gamma) \rightarrow \O_{\B_2} \rightarrow 0
\end{align*}
with $\O_W(K_W+\L_2-\Gamma)$,
we have
\begin{align}
\lefteqn{\dim H^0(W, \o_W^1(\log \B_2, \log C_3, \log C_3', \log \Gamma)(K_W+\L_2-\Gamma))} \nonumber\\
& & \le \dim H^0(W, \o_W^1(\log C_3, \log C_3', \log \Gamma)(K_W+\L_2-\Gamma))+1   \label{eq:chi22}
\end{align}
Note that $K_W+\L_2-\Gamma \equiv L-E_0-E_1'-E_2' \equiv C_1+E_1'+C_1'-E_2'$.
Applying the last statement of \lemref{lem:logsheaf} to $C_1, E_1'$ and $C_1'$,
we have
\begin{align}
\lefteqn{\dim H^0(W, \o_W^1(\log C_3, \log C_3', \log \Gamma)(K_W+\L_2-\Gamma))} \nonumber\\
& &=  \dim H^0(W, \o_W^1(\log C_3, \log C_3', \log \Gamma, \log C_1, \log E_1', \log C_1')(-E_2')) \label{eq:chi23}
\end{align}
The first Chern classes of $\Gamma, C_3, C_3', C_1, E_1', C_1'$ are linearly independent in $H^2(S, \mathbb{C})$.
By \cite{bidoublecover}*{Lemma~3.7}, it follows that
\begin{align}
\dim H^0(W, \o_W^1(\log C_3, \log C_3', \log \Gamma, \log C_1, \log E_1', \log C_1')(-E_2'))=0 \label{eq:chi24}
\end{align}
We thus obtain
$\dim H^2(S, \t_S)^{\mathrm{\chi_2}}=\dim H^0(W, \o_W^1(\log \Delta_2)(K_W+\L_2)) \le 2$
from \eqref{eq:chi21}-\eqref{eq:chi24}.\\

It remains to calculate $\dim H^0(W, \o_W^1(\log \Delta_3)(K_W+\L_3))$.
Note that $\Delta_3=\B_3 \cong \mathbb{P}^1$ and $\B_3(K_W+\L_3)=2$.
Tensoring the exact sequence
$$0 \rightarrow \o_W^1 \rightarrow \o_W^1(\log \B_3) \rightarrow \O_{\B_3} \rightarrow 0$$
with $\O_W(K_W+\L_3)$, we have
\begin{align}
\dim H^0(W, \o_W^1(\log \B_3))(K_W+\L_3)) \le \dim H^0(W, \o_W^1(K_W+\L_3))+3 \label{eq:chi31}
\end{align}
Note that $K_W+\L_3 \equiv (C_1+E_1'+C_1')+(C_2+E_2'+C_2')-E_3'$.
Applying the last statement of \lemref{lem:logsheaf} to $C_1, E_1', C_1', C_2, E_2'$ and $C_2'$, we have
\begin{align}
\lefteqn{\dim H^0(W, \o_W^1(K_W+\L_3))=} \nonumber \\
&& \dim H^0(W, \o_W^1(\log C_1, \log E_1', \log C_1', \log C_2, \log E_2', \log C_2')(-E_3'))\label{eq:chi32}
\end{align}
The Chern classes of $C_1, E_1', C_1', C_2, E_2'$ and $C_2'$ are linear independent.
By \cite{bidoublecover}*{Lemma~3.7}, it follows that
\begin{align}
H^0(W, \o_W^1(\log C_1, \log E_1', \log C_1', \log C_2, \log E_2', \log C_2')(-E_3'))=0  \label{eq:chi33}
\end{align}
We obtain
$\dim H^2(S, \t_S)^{\mathrm{\chi_3}}=h^0(W, \O_W(\log \Delta_3)(K_W+\L_3)) \le 3$ from \eqref{eq:chi31}-\eqref{eq:chi33}
and complete the proof.
\qed

\proof[The proof of \thmref{thm:moduli}]
From the construction of the surfaces in \exaref{exa:new},
it is clear that $\mathfrak{B}$ is irreducible, unirational and of dimension $3$
(see also \cite{chennew}*{Section~3}).

Let $S$ be a surface in \exaref{exa:new}.
Since $-\dim H^1(S,\t_{S})+\dim H^2(S,\t_{S})=2K_{S}^2-10\chi(S)=4$,
by \propref{prop:eigenspaces}, we have
$\dim H^2(S, \t_{S}) \le 7$ and $\dim H^1(S,\t_{S}) \le 3$.
Because $S$ is smooth and $K_S$ is ample,
the minimal model and the canonical model of $S$ coincide.
We denote by $\mathrm{B}(S)$ the base of the Kuranishi family of deformations of $S$
and by $[S]$ the corresponding point of the moduli space $\mathcal{M}^{\mathrm{can}}_{1,7}$.
Recall the fact that the germ of the complex space $(\mathcal{M}^{\mathrm{can}}_{1,7},[S])$
is analytically isomorphic to the quotient $\mathrm{B}(S)/\mathrm{Aut}(S)$.

We have the following inequalities
\begin{align*}
3 \ge \dim H^1(S,\t_S) &\ge \text{the dimension of $\mathrm{B}(S)$} \\
                  &= \text{the dimension of $(\mathcal{M}^{\mathrm{can}}_{1,7}, [S])$}\\
                  &\ge \text{the local dimension of $\mathfrak{B}$ at the point $[S]$} =3.
\end{align*}
Consequently all the equalities hold.
The second equality shows that $\mathrm{B}(S)$ is smooth.
Therefore $(\mathcal{M}^{\mathrm{can}}_{1,7},[S])$ is normal.
Since $\mathfrak{B}$ is irreducible,
the last equality shows that $\mathfrak{B}$ coincides with $\mathcal{M}^{\mathrm{can}}_{1,7}$
locally at $[S]$ for any $S$ in \exaref{exa:new}.
Hence $\mathfrak{B}$ is an open subset of $\mathcal{M}^{\mathrm{can}}_{1, 7}$ and $\mathfrak{B}$ is normal.

It remains to prove that $\mathfrak{B}$ is a closed subset of $\mathcal{M}^{\mathrm{can}}_{1, 7}$.
It suffices to prove the following statement.
{\itshape{Let $T$ be a smooth affine curve and $o \in T$.
Let $\mathcal{F} \colon \mathcal{S} \rightarrow T$ be a flat family of canonical models of surfaces
of general type with $K^2=7$ and $p_g=0$.
Set $\mathcal{S}_t:=\mathcal{F}^{-1}(t)$ for $t \in T$.
Assume that $\mathcal{S}_t$ is a surface in \exaref{exa:new} for $t \not =o$.
Then so is $\mathcal{S}_o$.}}

In fact, by construction, we have a $G$ ($\cong \mathbb{Z}_2^2$)-action on $\mathcal{S}\setminus \mathcal{S}_o$.
The $G$-action extends to $\mathcal{S}$ by \cite{moduli}*{Theorem~1.8}.
In particular, $\mathcal{S}_o$ admits a $G$-action.
This action lifts to the minimal model $\mathcal{S}_o'$ of $\mathcal{S}_o$.
By \propref{prop:basicK2=7}~(e), the canonical divisor of $\mathcal{S}_o'$ is ample and thus
$\mathcal{S}_o=\mathcal{S}_o'$.
It suffices to show that the pair $(\mathcal{S}_o, G)$ satisfies property (b) of \thmref{thm:classification}.

For $t \in T$,
we denote by $\mathcal{R}_{t}$ be the union of the divisorial parts of
the fixed loci of the three involutions $g_1, g_2$ and $g_3$.
Since $\mathcal{F} \colon \mathcal{S} \rightarrow T$ is a flat family with the $G$-action on each fiber,
$K_{\mathcal{S}_o} \mathcal{R}_{o} \ge K_{\mathcal{S}_t}\mathcal{R}_{t}$ for $t \not=o$.
From the construction in \exaref{exa:new}, we have $K_{\mathcal{S}_t} \mathcal{R}_{t}=5+5+3=13$ for $t \not =o$
and thus $K_{\mathcal{S}_o} \mathcal{R}_o \ge 13$.
We see that $(\mathcal{S}_o, G)$ satisfies the property (a) or (b) of \thmref{thm:classification}.

Set $\mathcal{Y}:=\mathcal{S}/G$.
Then $\mathcal{Y} \rightarrow T$ is a flat family of surfaces with only nodes.
For $t \not =o$, by construction, $\mathcal{Y}_t:=\mathcal{S}_t/G$ is a  $6$-nodal singular Del Pezzo surface.
It follows that $\mathcal{Y}_o:=\mathcal{S}_o/G$ has at least $6$ nodes.
Since the quotient of an Inoue surface is the $4$-nodal cubic surface,
by \thmref{thm:classification},
$(\mathcal{S}_o, G)$ satisfies property (b) and $\mathcal{S}_o$ is indeed a surface in \exaref{exa:new}.

Hence $\mathfrak{B}$ is a closed subset of $\mathcal{M}^{\mathrm{can}}_{1, 7}$
and we complete the proof of \thmref{thm:moduli}.
\qed.

\begin{rem}\label{rem:moduli}
From the proof of \thmref{thm:moduli}, we see that all the inequalities in \propref{prop:eigenspaces} hold,
$\dim H^1(S, \t_S)=7$ and $\dim H^2(S, \t_S)=3$.
We can also calculate the dimensions of the eigenspaces of $H^1(S, \t_S)$
(for the $G$-action):
$$\dim H^1(S, \t_S)^{\mathrm{inv}}=3\ \text{and}\ \dim H^2(S, \t_S)^{\mathrm{\chi_i}}=0\ \text{for}\ i=1, 2, 3.$$
It suffices to show that $\dim H^1(S, \t_S)^{\mathrm{inv}}=3$, i.e.,
$\dim H^1(W, \o_{W}^1(\log \Delta_1, \log \Delta_2, \log \Delta_3)(K_{W}))=3$ by \eqref{eq:eignesheaves}.
We have seen that
$\dim H^0(W, \o_{W}^1(\log \Delta_1, \log \Delta_2, \log \Delta_3)(K_W))=0$.
Note that $H^2(W, \o_{W}^1(\log \Delta_1, \log \Delta_2, \log \Delta_3)(K_W))$ is a
direct sum of $H^0(V, \t_V)$ by \eqref{eq:eignesheaves}.
Since $V$ is of general type, $H^0(V, \t_V)=0$.
It suffices to show that $\chi(\o_{W}^1(\log \Delta_1, \log \Delta_2, \log \Delta_3)(K_W))=-3$.
From the exact sequence \eqref{eq:newinv.log.seq},
we have
$$\chi(\o_{W}^1(\log \Delta_1, \log \Delta_2, \log \Delta_3)(K_{W})=\chi(\o_{W}^1(K_W))
+\chi(\oplus_{i=1}^3\O_{\Delta_i}(K_{W})).$$
The splitting principle and the Riemann-Roch theorem imply $\chi(\o_{W}^1(K_W))=-8$.
Each component of $\Delta_1, \Delta_2$ and $\Delta_3$ is a smooth rational curve (See \eqref{eq:newdata1}).
We obtain $\chi(\oplus_{i=1}^3\O_{\Delta_i}(K_W))=5$ by Riemann-Roch Theorem and
complete the proof.
\end{rem}

\paragraph{Acknowledgement.}Part of this work is from the author's Ph.~D thesis at the Peking University.
The author gratefully acknowledges the helpful suggestions of his former supervisor Jinxing~Cai.
The author thanks Xiaotao~Sun for all the support during the preparation of the paper.
The author is greatly indebted to Yifei~Chen, Wenfei~Liu, Jinsong~Xu and Lei~Zhang for many discussions.
The author wishes to thank Christopher Hacon for the invitation to the University of Utah
and for the hospitality.
This work is partially supported by the China Postdoctoral Science Foundation (Grant No.: 2013M541062).

\end{document}